\documentclass{article}
\usepackage{graphicx}
\usepackage[T1]{fontenc}
\usepackage[english]{babel}
\usepackage{graphicx}
\usepackage{amssymb, amsmath, amsfonts}
\usepackage{dsfont}
\usepackage{vmargin}
\usepackage{makeidx}

\everymath{\displaystyle}

\usepackage{hyperref}
\hypersetup{
	colorlinks=true,
	breaklinks=true,
	urlcolor=blue,
	linkcolor=red,
	bookmarksopen=true,
	pdftitle = Stochastic Cahn-Hilliard equation with singular nonlinearity and reflection,
	pdfauthor = Ludovic Gouden\`ege,
	pdfsubject = Stochastic Cahn-Hilliard equation with singular nonlinearity and reflection,
	pdfkeywords = {Cahn-Hilliard, stochastic partial differential equations, integration by parts formulae, reflection measures, invariant measures, singular nonlinearity}
}

\numberwithin{equation}{section} 

\newcommand{\R}{\mathbb{R}} 
\newcommand{\B}{\mathbf{B}} 
\newcommand{\N}{\mathbb{N}} 
\newcommand{\E}{\mathbb{E}} 
\newcommand{\PP}{\mathbb{P}} 
\newcommand{\F}{\mathcal{F}} 
\newcommand{\1}{\mathds{1}} 

\newcommand{\BD}{\begin{displaymath}}
\newcommand{\ED}{\end{displaymath}}
\newcommand{\BEA}{\begin{eqnarray}}
\newcommand{\EEA}{\end{eqnarray}}
\newcommand{\BEAS}{\begin{eqnarray*}}
\newcommand{\EEAS}{\end{eqnarray*}}
\newcommand{\BE}{\begin{equation}}
\newcommand{\EE}{\end{equation}}
\newcommand{\BES}{\begin{equation*}}
\newcommand{\EES}{\end{equation*}}

\newcommand{\Proof}{\noindent {\bf{\underline{Proof} : }}}
\newcommand{\EProof}{\begin{flushright}$\Box$\end{flushright}}

\newcommand{\Step[1]}{\indent \textbf{Step {#1}.}}
\newcommand{\EStep[1]}{\begin{flushright}$\Box$\end{flushright}}
\newcommand{\Case[1]}{\textbf{Case {#1}}}

\newcommand{\FAN}{\quad \text{ for all }n\geq 1}

\newcommand{\FAX}{\quad \text{ for all }x \in L^2}

\newcommand{\fln}{f_{\ln}}
\newcommand{\fnln}{f^n_{\ln}}
\newcommand{\Fln}{F_{\ln}}
\newcommand{\Fnln}{F^n_{\ln}}
\newcommand{\Foneln}{F^1_{\ln}}

\newtheorem{Def}{Definition}[section]
\newtheorem{Coro}{Corollary}[section]
\newtheorem{Le}{Lemma}[section]
\newtheorem{Th}{Theorem}[section]
\newtheorem{Prop}{Proposition}[section]
\newtheorem{Rem}{Remark}[section]

\title{Stochastic Cahn-Hilliard equation with singular nonlinearity and reflection}
\author{Ludovic Gouden\`ege\\
{\small ENS Cachan - Antenne de Bretagne, Campus de Ker-Lann, 35170 Bruz, France}\\
{\small ludovic.goudenege@bretagne.ens-cachan.fr}}
\date{}

\begin{document}

\maketitle
\footnotetext{\!\!\!\!\!\!\!\!\!\!\!AMS 2000 subject classifications. 60H15, 60H07, 37L40.\\
{\em Key words and phrases} : Cahn-Hilliard, stochastic partial differential equations, integration by parts formulae, reflection measures, invariant measures, singular nonlinearity.}
\tableofcontents
\newpage

\section*{R\'esum\'e}
On consid\`ere une \'equation aux d\'eriv\'ees partielles stochastique poss\'edant une non-lin\'earit\'e de type logarithmique (ou une puissance n\'egative), avec une reflexion en z\'ero sous la contrainte de conservation de masse. L'\'equation, dirig\'ee par un bruit blanc en espace et en temps, contient un double Laplacien. L'absence de principe de maximum pour le double Laplacien pose des difficult\'es pour l'utilisation d'une m\'ethode classique de p\'enalisation, pour laquelle une importante propri\'et\'e de monotonie est utilis\'ee. Etant inspir\'e par les travaux de Debussche et Zambotti, on emploie une m\'ethode bas\'ee sur les \'equations en dimension infinie, utilisant l'approximation par des \'equations reguli\`eres et la convegence des semi-groupes de transition li\'es aux \'equations r\'egularis\'ees. On d\'emontre l'existence et l'unicit\'e de solutions pour des donn\'ees initiales positives, et on donne plusieurs resultats sur les mesures invariantes et les mesures de r\'eflexion.
\section*{Abstract}
We consider a stochastic partial differential equation with logarithmic (or negative power) nonlinearity, with one reflection at 0 and with a constraint of conservation of the space average. The equation, driven by the derivative in space of a space-time white noise, contains a bi-Laplacian in the drift. The lack of the maximum principle for the bi-Laplacian generates difficulties for the classical penalization method, which uses a crucial monotonicity property. Being inspired by the works of Debussche and Zambotti, we use a method based on infinite dimensional equations, approximation by regular equations and convergence of the approximated semi-group. We obtain existence and uniqueness of solution for nonnegative intial conditions, results on the invariant measures, and on the reflection measures.
\section*{Introduction and main results}
The Cahn-Hilliard-Cook equation is a model to describe phase separation in a binary alloy (see \cite{C}, \cite{CH1} and \cite{CH2}) in the presence of thermal fluctuations (see \cite{COOK} and \cite{LANGER}). It takes the form:
\BE\label{Eq:0.1}
\left\{\begin{array}{ll}                                                                                                                                                                                                                                                                                                       
\partial_{t} u = - \frac{1}{2}\Delta\left(\Delta u - \psi(u)\right)+\dot{\xi},&\text{ on } \Omega\subset \R^n,\\
\\
\nabla u \cdot \nu = 0 = \nabla (\Delta u) \cdot \nu, &\text{ on } \partial\Omega,\\                                                                                                                                                                                                      
\end{array}\right.
\EE                                                                                                                                                                                                                                                                                                                     
where $t$ denotes the time variable and $\Delta$ is the Laplace operator. Also $u\in [-1,1]$ represents the ratio between the two species and the noise term $\dot{\xi}$ accounts for the thermal fluctuations. The nonlinear term $\psi$ has the double-logarithmic form:                      
\BE\label{Eq:0.2}                                                                                                                                                                                                                                                                                                       
\psi : u \mapsto \ln\left(\frac{1+u}{1-u}\right) - \kappa u.                                                                                                                                                                                                                                                                       
\EE                                                                                                                                                                                                                                                                                                                     
The deterministic equation has been extensively studied first in the case where $\psi$ is replaced by a polynomial function (see \cite{CH1}, \cite{LANGER} and \cite{MR763473}) and then for non smooth $\psi$ (see \cite{MR1123143} and \cite{MR1327930}). Furthermore, this model has been used successfully for describing phase separation phenomena, see for example the survey \cite{MR1657208}, and the references therein, or others recent results on spinodal decomposition and nucleation in \cite{MR1232163, MR2342011, MR1214868, MR1637817, MR1753703, MR1712442, MR1763320, MR2048517}. In the polynomial case, the concentration $u$ is not constricted to remain between $-1$ and $1$ and the logarithmic nonlinearity might seem preferable.\\
Up to our knowledge, only the polynomial nonlinearity has been studied in the stochastic case (see \cite{BLMAWA1, BLMAWA2, MR1867082, MR1897915, MR1359472, MR1111627}). This article is a step toward the mathematical comprehension of the full model with double-logarithmic term and noise. We consider the one dimensional case and consider a nonlinear term with only one singularity.
Clearly, due to the noise, such an equation cannot have a solution, and a reflection measure should be added to the equation. Thus the right stochastic equation to study is:
\BE\label{Eq:0.10}                                                                                                                                                                                                                                                                                                      
\left\{\begin{array}{ll}                                                                                                                                                                                                                                                                                                
\displaystyle{\partial_{t} X = - \frac{1}{2}\Delta\Big(\Delta X + f(X) + \eta\Big) + \partial_{\theta}\dot{W}},&\text{ with } \theta \in [0,1]=\Omega,\\                                                                                                                                                                                                     
&\\                                                                                                                                                                                                                                                                                                                     
\nabla X \cdot \nu = 0 = \nabla (\Delta X) \cdot \nu, &\text{ on } \partial\Omega,
\end{array}\right.                                                                                                                                                                                                                                                                                                      
\EE
where $f$ is defined below, and where the measure is subject to the contact condition almost surely:
\BE\label{Eq:0.9}                                                                                                                                                                                                                                                                                           
\int X d\eta = 0.                                                                                                                                                                                                                                                                                              
\EE   
Stochastic partial differential equations with reflection can model the described problem or the evolution of random interfaces near a hard wall (see \cite{MR1835843} and \cite{MR2128236}). For other results on fluctuations of random interfaces, see \cite{MR1872740}. For a detailled study of the contact set $\{(t,\theta) : X(t,\theta)=0\}$ and of the reflection measure $\eta$, see \cite{MR2257651}, \cite{MR1921014} and \cite{MR1959795}. The equation \eqref{Eq:0.10} has been studied when no nonlinear term is taken into account in \cite{MR2349572}. In this paper, the authors have introduced various techniques needed to overcome the lack of comparison principle for fourth order equations. Indeed, the case of a second order equation was studied in \cite{MR1172940} where an extensive use of monotonicity is used, as well as in all the articles treating with the second order case.\\
This article is in the spirit of \cite{MR1959795} where a nonlinear term is taken into account for the second order equation. We study existence and uniqueness of solution for equation \eqref{Eq:0.10} with $f$ of the form:
\BE\label{Eq:0.11}                                                                                                                                                                                                                                                                                                      
f(x):=\fln(x):=\left\{\begin{array}{lc}                                                                                                                                                                                                                                                                                 
-\ln x,& \text{ for all } x >0\\                                                                                                                                                                                                                                                                                        
+\infty,& \text{ for all } x \leq 0,\\                                                                                                                                                                                                                                                                                  
\end{array}\right.                                                                                                                                                                                                                                                                                                      
\EE                                                                                                                                                                                                                                                                                                                     
or for $\alpha>0$:                                                                                                                                                                                                                                                                                                    
\BE\label{Eq:0.12}                                                                                                                                                                                                                                                                                                      
f(x):=f_{\alpha}(x):=\left\{\begin{array}{lc}                                                                                                                                                                                                                                         
x^{-\alpha},& \text{ for all } x >0\\                                                                                                                                                                                                                                                                                   
+\infty,& \text{ for all } x \leq 0.\\                                                                                                                                                                                                                                                                                  
\end{array}\right.                                                                                                                                                                                                                                                                                                      
\EE
Moreover we characterize the case when the measure $\eta$ vanishes. Our method mixes ideas from \cite{MR2349572} and \cite{MR1959795}. Additional difficulties are overcome, the main one being to understand how to deal with the nonlinear term. Again in \cite{MR1959795}, this term is not difficult to consider thanks to monotonicy arguments.\\
Our main results state that equations \eqref{Eq:0.10}, \eqref{Eq:0.9} together with an initial condition have an unique solution (see \ref{Th:3.1} and \ref{Th:Existence}). It is constructed thanks to the gradient structure of \eqref{Eq:0.10} and Strong Feller property. Furthermore, we prove that the measure $\eta$ vanishes only for $f$ described in \eqref{Eq:0.12} with $\alpha \geq 3$ (see \ref{Th:5.4}).

\newpage
\section{Preliminaries}\label{S:1}
\subsection{Notation}
We denote by $\langle\cdot,\cdot\rangle$ the scalar product in $L^2(0,1)$: 
\BD
\text{for all } h,k \in L^2(0,1) \qquad \langle h,k\rangle = \int_{0}^1h(\theta) k(\theta)d\theta.
\ED
We denote by $A$ the realization in $L^2(0,1)$ of the Laplace operator with Neumann boundary condition, i.e.: 
\BD
D(A) = \text{ Domain of }A  = \{ h \in W^{2,2}(0,1) : h'(0) = h'(1) = 0\}
\ED
where the space $W^{2,2}(0,1)$ is the classical Sobolev space. Below we use the notation $W^{n,p}$  and $||.||_{W^{n,p}}$ to denote the Sobolev space $W^{n,p}(0,1)$ and its associated norm.
Remark that $A$ is self-adjoint on $L^2(0,1)$ and we have a complete orthonormal system of eigenvectors $(e_i)_{i \in \N}$  in $L^2(0,1)$ associated to the eigenvalues $\lambda_{i}:=(-(i\pi)^2)_{i\in\N}$ where we define:
\BD
e_0(\theta) =1,\ e_i(\theta) = \sqrt{2}\cos(i\pi\theta),\text{ for all } i \in \N^*,\text{ for all }\theta \in [0,1].
\ED
We denote by $\bar{h}$ the mean of $h \in L^2(0,1)$: 
\BD
\bar{h} = \int_{0}^1h(\theta)d\theta = \langle h,e_0\rangle.
\ED
Then we define for all $c \in \R$ : 
\BD
L^2_c = \{h\in L^2(0,1) :\bar{h} = c\},
\ED
and $L^2=L^2(0,1)$. We remark that $(-A)^{-1} : L^2_0\rightarrow L^2_0$ is well defined. We denote by $Q$ this operator. We can extend the definition of $Q$ to $L^2(0,1)$ (we denote this operator $\bar{Q}$) by the formula: 
\BD
\bar{Q}h = Q(h-\bar{h})+\bar{h}, \text{ for all } h \in L^2(0,1)
\ED
For $\gamma \in\R$, we define $(-A)^{\gamma}$ by setting 
\BD
(-A)^{\gamma} h = \sum_{i=1}^{+\infty} (-\lambda_{i})^{\gamma} h_{i}e_{i}, \text{ when } h = \sum_{i=0}^{+\infty} h_{i}e_{i}.
\ED
The domain of $(-A)^{\gamma/2}$ is
\BD
V_{\gamma}:=D((-A)^{\gamma/2}) = \Big\{h = \sum_{i=0}^{+\infty} h_{i}e_{i} :\sum_{i=1}^{+\infty} (-\lambda_{i})^{\gamma} h_{i}^2<+\infty\Big\}.
\ED
It is endowed with the seminorm 
\BD
|h|_{\gamma} =  \left(\sum_{i=1}^{+\infty} (-\lambda_{i})^{\gamma} h_{i}^2\right)^{1/2},
\ED
and with the norm
\BD
\|h\|_{\gamma} =  \big(|h|_{\gamma}^2 + \bar{h}^2\big)^{1/2},
\ED
associated to the scalar product defined for all $h,k \in V_{\gamma}$ by $( h,k)_{\gamma}$.\\
For $\gamma=-1$, $V_{-1}=D((-A)^{-1/2})$ is the completion of the space of functions $h \in L^2$ such that
\BEAS
\langle \bar{Q}h,h\rangle = \langle Q(h-\bar{h})+\bar{h},h\rangle &=& \langle (-A)^{-1}(h-\bar{h}),h-\bar{h}\rangle + \bar{h}^2\\
&=& \langle (-A)^{-1/2}(h-\bar{h}),(-A)^{-1/2}(h-\bar{h})\rangle\\
&=&  |h|_{-1}^2 + \bar{h}^2 < +\infty.
\EEAS
To lighten notations, we set $(\cdot,\cdot) := (\cdot,\cdot)_{-1}$ for the inner product of $V_{-1}$. The average plays an important role and we often work with functions with a fixed average $c\in \R$. We define $H_c=\{h\in H , \bar{h}=c\}$ for all $c\in \R$. We set
\BD
D(B) = W_0^{1,2}(0,1), B = \frac{\partial}{\partial \theta}, D(B^*) = W^{1,2}(0,1) \text{ and } B^*=-\frac{\partial}{\partial \theta}.
\ED
We remark that $BB^*=-A$. Finally, we denote by $\Pi$ the orthogonal projector of $V_{-1}$ onto $H_0$. We have: 
\BES
\begin{array}{rcl}
\Pi : V_{-1} &\rightarrow& H_0\\
h & \mapsto & h-\bar{h}.
\end{array}
\EES
Notice that $\Pi$ is also an orthogonal projector of $L^2$ onto $L^2_0$. Moreover:
\BE\label{Eq:1.1}
-A\bar{Q}h = \Pi h, \text{ for all } h \in L^2(0,1).
\EE
We denote by $\mathcal{B}_b(H_{c})$ the space of all Borel bounded functions and $\mathcal{C}_b(H_{c})$ the space of continous bounded functions. We set $O_{s,t} := [s,t] \times [0,1]$ for $s,t \in [0,T]$ with $s<t$ and $T>0$, and $O_{t} = O_{0,t}$ for $0\leq t\leq T$. Given a measure $\zeta$ on $O_{s,t}$ and a continuous function $v$ on $O_{s,t}$, we set
\BD
\big\langle v,\zeta\big\rangle_{O_{s,t}} := \int_{O_{s,t}} v d\zeta.
\ED
In order to solve the equation \eqref{Eq:0.10}, we use a Lipschitz approximation of this equation. We denote by $\{f^n\}_{n \in \N}$ the sequence of Lipschitz functions which converges to the function $f$ on $(0,+\infty)$, defined for $n \in \N$ by: 
\BD
f^n(x):=f(x^++1/n),\text{ for all } x \in \R.
\ED
When $f=\fln$ is the logarithmic function \eqref{Eq:0.11}, we use the following positive antiderivative of $-f^n=-\fnln$
\BD
F^n(x)=\Fnln(x):=(x+1/n) \ln(x^++1/n)-x^++1-1/n,\text{ for all } x \in \R,
\ED
and the following positive antiderivative of $-f=-\fln$ defined only on $\R^+$ by:
\BD
F(x)=\Fln(x):=x \ln(x)-x+1,\text{ for all }x \in \R^+.
\ED
When $f=f_{\alpha}$ is the negative $\alpha$-power function \eqref{Eq:0.12} with $\alpha \neq 1$, we use the following antiderivative of $-f^n=-f^n_{\alpha}$
\BD
F^n(x)=F^n_{\alpha}(x):=\frac{(x^++1/n)^{1-\alpha}}{\alpha-1} + n^{\alpha}x^-,\text{ for all } x \in \R,
\ED
and the following antiderivative of $-f=-f_{\alpha}$ defined only on $\R^+$ by:
\BD
F(x)=F_{\alpha}(x):=\frac{x^{1-\alpha}}{\alpha-1},\text{ for all } x \in \R^+.
\ED
Finally when $\alpha=1$, we use the following antiderivative of $-f^n=-f^n_{\alpha}$
\BD
F^n(x)=F^n_{\alpha}(x):=-\ln(x^++1/n) + n x^-,\text{ for all } x \in \R,
\ED
and the following antiderivative of $-f=-f_{\alpha}$ defined only on $\R^+$ by:
\BD
F(x)=F_{\alpha}(x):=-\ln x,\text{ for all } x \in \R^+.
\ED
We use the notation $f,f^n,F,F^n$ when the result holds both for $\fln$ and $f_{\alpha}$. Otherwise we use $\fln,\fnln,\Fln,\Fnln$ or $f_{\alpha},f^n_{\alpha},F_{\alpha},F^n_{\alpha}$.\\
With these notations, we rewrite \eqref{Eq:0.10} in the abstract form:
\BE\label{Eq:1.2}
\left\{\begin{array}{l}
dX=-\frac{1}{2}(A^2X+Af(X))dt + BdW,\\
\\
\langle X ,\eta\rangle_{O_{T}} = 0,\\
\\
X(0,x)=x \text{ for } x \in V_{-1}.
\end{array}\right.
\EE
Finally, in all the article, C denotes a constant which may depend on $T$ and $\alpha$ and its value may change from one line to another.

\subsection{The linear equation}\label{S:2}
The linear equation is given by
\BES
\left\{\begin{array}{l}
dZ(t,x) = -\frac{1}{2} A^2Z(t,x)dt + BdW, \quad \text{ for all } t \in [0,T],\\
\\
Z(0,x) = x.
\end{array}\right.
\EES
where $x \in V_{-1}$. We have
\BES
Z(t,x) = e^{-tA^2/2}x + \int_0^te^{-(t-s)A^2/2}BdW_s.
\EES
As easily seen this process is in $\mathcal{C}([0,+\infty[;L^2(0,1))$ (see \cite{MR1207136}). In particular, the mean of $Z$ is constant and the law of the process $Z(t,x)$ is the Gaussian measure:
\BD
Z(t,x)\sim\mathcal{N}\big(e^{-tA^2/2}x,Q_t\big),
\ED
where
\BD
Q_t = \int_0^te^{-sA^2/2}BB^*e^{-sA^2/2}ds = (-A)^{-1}(I - e^{-tA^2}).
\ED
\EProof
If we let $t\rightarrow +\infty$, the law of $Z(t,x)$ converges to the Gaussian measure on $L^2$:
\BD
\mu_c := \mathcal{N}(ce_0,Q), \text{ where } c=\bar{x}.
\ED
Notice that the kernel of $Q$ is $\{te_0,t\in\R\}$ and $\mu_c$ is concentrated on $L^2_c$. It is important to remark that the measure $\mu_c$ is linked to the Brownian motion. Indeed, let $(\B_{\theta})_{\theta\in[0,1]}$ be a Brownian motion, then the law of $Y_{c}(\theta) = \B(\theta) - \overline{\B} + c$ is $\mu_c$ (see \cite{MR2349572}).

\subsection{Lipschitz Approximation}\label{S:3}
For $n \in \N$, we study for the following Lipschitz approximation of \eqref{Eq:1.2} with an initial condition $x \in V_{-1}$: 
\BE\label{Eq:2.1}
\left\{\begin{array}{l}
dX^n+\frac{1}{2}(A^2X^n+Af^n(X^n))dt = BdW,\\
\\
X^n(0,x)=x.
\end{array}\right.
\EE
We prove existence and uniqueness of solution in a suitable space for the equation \eqref{Eq:2.1}. We then follow standard arguments to show existence and uniqueness of an invariant measure for the equation \eqref{Eq:2.1} with fixed $n \in \N$, and the strong Feller property of the semigroup. First we have to define the definition of a weak solution to \eqref{Eq:2.1}.\\
We say $X^n$ is a mild solution of \eqref{Eq:2.1} if it is satisfied for all $t\geq 0$:
\BE\label{Eq:2.5}
X^n(t,x)=Z(t,x)-\int_0^t Ae^{-(t-s)A^2/2}f^n(X^n(s,x))ds.
\EE
\begin{Le}\label{Le:2.1}
Fix $n\in\N$, $0<\varepsilon<2/3$ and $p=4(1-\varepsilon)$. For all $x \in L^2(0,1)$ there exists a unique adapted process $X^n \in \mathcal{C}([0,T];V_{-1}) \cap L^{p}([0,T];L^2(0,1))$ solution of equation \eqref{Eq:2.5}. Moreover for all $t\geq0$:
\BE\label{Eq:2.6}
\langle X^n(t,x),e_{0}\rangle = \langle x, e_{0}\rangle.
\EE
\end{Le}
\Proof The proof is classical and left to the reader. It is based on the following inequalities
\BEA
||(-A)^{1/2} e^{-tA^2/2}h||_{0} &\leq& C ||h||_{0}t^{-1/4}, \quad t>0,\quad h \in L^2\\
||A e^{-tA^2/2}h||_{0} &\leq& C ||h||_{0}t^{-1/2}, \quad t>0,\quad h \in L^2\\
||e^{-tA^2/2}h||_{0} &\leq& C |h|_{-1}t^{-1/4}, \quad t>0,\quad h \in L^2.
\EEA
\EProof
It is also standard to prove 
\begin{Le}\label{Le:2.2}
For $n \in \N$ and $c \in \R$, for all $t>0$:
\BE\label{Eq:2.12}
|X^n(t,x)-X^n(t,y)|_{-1} \leq \exp(-t\pi^4/2) |x-y|_{-1},\quad \text{ for all } x,y \in L^2_c.
\EE
\end{Le}
\Proof
We consider for $N\in \N$ and $x,y \in L^2_c$ the process :
\BD
S^{N}(t,x,y) = \sum_{i=0}^N \langle X^n(t,x) - X^n(t,y),e_{i} \rangle e_{i},\quad \text{ for all } t \geq 0.
\ED
then $t \mapsto S^N(t,x,y)$ is $\mathcal{C}^1$ with values in a $(N+1)$-dimensional subspace of $D(A)$ such that $\Pi S^N = S^N$. Indeed 
\BEAS
\Pi S^N(t,x,y) &=& \sum_{i=1}^N \langle X^n(t,x) - X^n(t,y),e_{i} \rangle \Pi e_{i} + \langle X^n(t,x) - X^n(t,y),e_{0} \rangle \Pi e_{0}\\
&=& \sum_{i=1}^N \langle X^n(t,x) - X^n(t,y),e_{i} \rangle e_{i}\\
&=& S^N - \langle X^n(t,x) - X^n(t,y),e_{0} \rangle e_{0}. 
\EEAS
And by \eqref{Eq:2.6} and since $x$ and $y$ are in $L^2_c$
\BES
\langle X^n(t,x) - X^n(t,y),e_{0} \rangle  = \langle x - y,e_{0} \rangle = 0.
\EES
By \eqref{Eq:1.1} we have $-A\bar{Q}h= \Pi h$ for all $h \in L^2(0,1)$, then $-A\bar{Q}S^N = S^N$. Using the spectral behavior of $A$ given in section \ref{S:1}, we have the following computation: 
\BEAS
\frac{d}{dt}|S^N(t,x,y)|_{-1}^2 &=& \frac{d}{dt}\langle \bar{Q}S^N(t,x,y),S^N(t,x,y)\rangle\\
&=&2\langle \frac{d}{dt}S^N(t,x,y),\bar{Q}S^N(t,x,y)\rangle\\&&\\
&=&\langle -A^2S^N(t,x,y),\bar{Q}S^N(t,x,y)\rangle\\
&&+ \langle -A(f^n(X^n(t,x))-f^n(X^n(t,y))),\bar{Q}S^N(t,x,y)\rangle\\&&\\
&=&\langle AS^N(t,x,y),S^N(t,x,y)\rangle\\
&&+ \langle f^n(X^n(t,x))-f^n(X^n(t,y)),S^N(t,x,y)\rangle\\&&\\
&\leq& -\pi^4|S^N(t,x,y)|_{-1}^2 \\
&&+ \langle f^n(X^n(t,x))-f^n(X^n(t,y)),S^N(t,x,y)\rangle
\EEAS
This differential inequality implies : 
\BES
|S^N(t,x,y)|_{-1}^2 \leq e^{-t\pi^4}|x-y|_{-1}^2
+ \int_{0}^t e^{-(t-s)\pi^4}\langle f^n(X^n(s,x))-f^n(X^n(s,y)),S^N(s,x,y)\rangle ds.
\EES
Moreover by letting $N \rightarrow +\infty$ we have $|S^N(t,x,y)|_{-1}^2\rightarrow |X^n(t,x)-X^n(t,y)|_{-1}^2$, and since $f^n$ is monotone non-increasing we obtain 
\begin{multline*}
\langle f^n(X^n(s,x))-f^n(X^n(s,y)),S^N(s,x,y)\rangle\\
\mathop{\longrightarrow}_{N\rightarrow+\infty} \langle f^n(X^n(s,x))-f^n(X^n(s,y)),X^n(s,x)-X^n(s,y)\rangle \leq 0.
\end{multline*}
Then the limit of the integral is nonpositive, and we obtain the expected inequality \eqref{Eq:2.12}.
\EProof
It is classical that $X^n \in \mathcal{C}([0,T];V_{-1}) \cap L^{p}([0,T];L^2(0,1))$ satisfies \eqref{Eq:2.5} if and only if it is a weak solution of \eqref{Eq:2.1} in the sense
\begin {Def}\label{De:1}
For $n\in\N$, $0<\varepsilon<2/3$ and $p=4(1-\varepsilon)$, let $x \in \mathcal{C}([0,1],\R^+)$ with $\overline{x} > 0$. We say that $(X^n(t,x))_{t\in[0,T]}$, defined on a stochastic basis linked to $(W(t))_{t\in [0,T]}$, is a solution to \eqref{Eq:2.1} on $[0,T]$ if :
\begin{enumerate}
\item[(a)] almost surely $X^n(\cdot,x) \in \mathcal{C}([0,T];V_{-1})\cap  L^{p}([0,T];L^2(0,1))$,
\item[(b)] for all $h \in D(A^2)$ and for all $0 \leq t \leq T$ :
\BES
\begin{array}{lll}
\langle X^n(t,x),h\rangle&=&\displaystyle\langle x,h\rangle - \int_{0}^t\langle X^n(s,x),A^2h\rangle ds\\
&&\displaystyle -\int_{0}^t\langle Ah,f^n(X^n(s,x))\rangle ds - \int_{0}^{t}\langle Bh, dW\rangle.
\end{array}
\EES
\end{enumerate}
\end{Def}
We now describe an important property of equation \eqref{Eq:2.1}. It can be described as a gradient system in $V_{-1}$ with a convex potential, and can be rewritten as:
\BE\label{Eq:2.2}
\left\{\begin{array}{l}
dX^n-\frac{1}{2}A(-AX^n+\nabla U^n(X^n))dt = BdW,\\
\\
X^n(0,x)= x \in L^2(0,1),
\end{array}\right.
\EE
where $\nabla$ denotes the gradient in the Hilbert space $L^2(0,1)$, and : 
\BE
U^n(x):= \int_{0}^{1}F^n(x(\theta))d\theta, \quad x \in L^2(0,1).
\EE
Notice that $\nabla U^n(x) = -f^n(x)$ which is dissipative, then $U^n$ is a convex potential. Finally, we define the probabilty measure on $L^2_{c}$:
\BE\label{Eq:2.10}
\nu_c^n(dx)=\frac{1}{Z_c^n}\exp(-U^n(x))\mu_c(dx),
\EE
where $Z_c^n$ is a normalization constant. By Lemma \ref{Le:2.2}, we easily obtain that the equation $\eqref{Eq:2.1}$ in $H_{c}$ has a unique ergodic invariant measure and it is not difficult to prove that this measure is precisely $\nu_{c}^n$. Since the potential $U^n$ is convex, we can prove that the transition semigroup is strong Feller. Let $(P^{n,c}_t)_{n \in \N}$ be the sequence of transition semigroup for an initial condition in $H_{c}$ such that
\BES
P^{n,c}_t\phi(x) = \E[\phi(X^{n,c}(t,x)], \quad\text{ for all } t \geq 0, x \in H_{c}, \phi \in \mathcal{B}_b(H_{c})\text{ and } n \in \N^*,
\EES
where $X^{n,c}(t,x)$ is the solution of the equation \eqref{Eq:2.2}.
\begin{Prop}\label{Prop:2.1} For abitrary $T>0$, there exists a constant $C_{T} > 0$ such that for all $\phi \in \mathcal{B}_{b}(H_{c})$, for all $n\in\N$ and for all $t \in [0,T]$: 
\BE\label{Eq:2.14}
| P^{n,c}_{t}\phi(x) - P^{n,c}_{t}\phi(y)| \leq \frac{\sqrt{C_{T}}}{\sqrt{t}}\|\phi\|_\infty\|x-y\|_{-1}, \quad\text{for all } x,y \in H_{c}.
\EE
\end{Prop}
\Proof
We now consider the following process :
\BES
\begin{array}{rcl}
H_{0} & \rightarrow & H_{0}\\
x & \mapsto & \mathcal{X}^{n,c}(t,x) = X^n(t,x + ce_{0}) - ce_{0} 
\end{array}
\EES
which solves the following equation : 
\BE\label{Eq:2.13}
\left\{\begin{array}{l}
\displaystyle{d\mathcal{X}^{n,c} - \frac{1}{2}A(-A\mathcal{X}^{n,c}+\nabla U^n(ce_{0}+\mathcal{X}^{n,c}))dt = BdW},\\
\displaystyle{\mathcal{X}^{n,c}(0,x)=x \in H_{0}}.
\end{array}\right.
\EE
This equation describes a gradient system in $H_{0}$ with non-degenerate noise and with a convex potential. We fix $c>0$ and $n\in\N$, and denote $\mathcal{P}^{n,c}_{t}$ the markov transition semigroup defined by :
\BES
\mathcal{P}^{n,c}_{t}\psi(x) = \E[\psi(\mathcal{X}^{n,c}(t,x))],\quad \text{for all } t \geq 0, \text{for all } x \in H_{0}, \text{for all } \psi \in \B_b(H_{0}).
\EES
For all $c\in\R$, for all $x\in H_{c}$ and $\psi \in \mathcal{B}_b(H_{0})$, if we set 
\BD
\begin{array}{rcl}
\phi : H_{c} &\rightarrow &\R\\
u &\mapsto& \psi (u - ce_{0}),
\end{array}
\ED
we have the following equality:
\BE\label{Eq:2.13b}
\mathcal{P}^{n,c}_{t}\psi(x-ce_{0}) = P^{n,c}_{t}\phi(x).
\EE
Then for all $c\in\R$, for all $x\in H_{c}$ and $\psi \in \mathcal{B}_b(H_{0})$, the following Bismut-Elworthy formula holds\footnote{As usual some computations below are formal and would be difficult to justify rigourously in our infinite dimensionnal setting. However the final result is easy to justify by Galerkin approximation (see \cite{MR1359472}, section 3.2).}:
\BE\label{Eq:2.15}
D\mathcal{P}_{t}^{n,c}\psi(x-ce_{0})\cdot h = \frac{1}{t}\E\left[\psi(\mathcal{X}^{n,c}(t,x-ce_{0}))\int_{0}^t \langle(D\mathcal{X}^{n,c}(s,x-ce_{0})\cdot h),dW\rangle\right].
\EE
Then by \eqref{Eq:2.13} and \eqref{Eq:2.15},
\BE\label{Eq:2.16}
| DP_{t}^{n,c}\phi(x)\cdot h|^2 \leq \frac{1}{t^2}\|\phi\|_{\infty}^2\E\left[\int_0^t\|DX^{n,c}(s,x)\cdot h\|_{-1}^2ds\right].
\EE
Let $x$ and $y$ be arbitrary elements in $H_{c}$, then by the mean value theorem, for $\sigma(y) \in [0,T]$
\BE\label{Eq:2.20}
P_{t}^{n,c}\phi(x) - P_{t}^{n,c}\phi(y) = DP_{t}^{n,c}\phi(x+\sigma(y)(x-y))\cdot (x-y).
\EE
We use an estimate on $\|DX^{n,c}(s,x+\sigma(y)(x-y))\cdot (x-y)\|_{-1}^2$, \eqref{Eq:2.16} and \eqref{Eq:2.20}, and we have the expected result for all $x$, $y \in H_{c}$
\BES
|P_{t}^{n,c}\phi(x) - P_{t}^{n,c}\phi(y)| \leq \frac{\sqrt{C_{T}}}{\sqrt{t}}\|\phi\|_{\infty}\|x-y\|_{H}.
\EES
\EProof

\newpage
\section{Solutions of equation with a reflection measure}\label{S:4}
For all $n \in \N$ we have a unique solution $X^{n}$ of \eqref{Eq:2.1}. We want to know if these solutions converge to a solution of the equation \eqref{Eq:0.10}. First we describe the definition of a weak solution for \eqref{Eq:0.10} :
\begin {Def}\label{De:2}
Let $x \in \mathcal{C}([0,1],\R^+)$ and $\overline{x} > 0$. We say that $\left(\left(X(t,x)\right)_{t\in[0,T]},\eta,W\right)$, defined on a filtered complete probability space $\left(\Omega, \PP, \F, (\F_t)_{t\in [0,T]}\right)$, is a weak solution to \eqref{Eq:0.10} on $[0,T]$ for the initial condition $x$ if:
\begin{enumerate}
\item[(a)] a.s. $X \in \mathcal{C}(]0,T]\times[0,1];\R^+) \cap \mathcal{C}([0,T];V_{-1})$ and $X(0,x) = x$,
\item[(b)] a.s. $\eta$ is a positive measure on $(0,T]\times [0,1]$, such that $\eta(O_{\delta,T}) <+\infty$ for all $\delta \in (0,T]$,
\item[(c)] $W$ is a cylindrical Wiener process on $L^2(0,1)$,
\item[(d)] the process $\left(X(\cdot,x),W\right)$ is $(\F_t)$-adapted,
\item[(e)] a.s. $f(X(\cdot,x)) \in L^1(O_{T})$,
\item[(f)] for all $h \in D(A^2)$ and for all $0 < \delta \leq t \leq T$ :
\BD
\begin{array}{rcl}
\langle X(t,x),h\rangle &=& \langle X(\delta,x),h\rangle - \int_{\delta}^t\langle X(s,x),A^2h\rangle ds-\int_{\delta}^t\langle Ah(\theta),f(X(s,x))\rangle ds\\
&&- \big\langle Ah,\eta\big\rangle_{O_{\delta,t}} - \int_{\delta}^t\langle Bh,dW\rangle, \quad a.s.,
\end{array}
\ED
\item[(g)] a.s. the contact property holds : $supp(\eta) \subset \{(t,\theta)\in O_{T} / X(t,x)(\theta)=0\}$, that is, 
\BD
\big\langle X,\eta\big\rangle_{O_{T}}= 0.
\ED
\end{enumerate}
Finally, a weak solution $(X,\eta,W)$ is a strong solution if the process $t \mapsto X(t,x)$ is adapted to the filtration $t \mapsto \sigma(W(s,.),s\in [0,t])$
\end{Def}

\begin{Rem} In (f), the only term where we use the function $f$ is well defined. Indeed, by (e) we have $f(X(\cdot,x)) \in L^1(O_{T})$ and by Sobolev embedding $Ah \in D(A) \subset L^\infty(O_{T})$. Hence the notation $\langle \cdot, \cdot \rangle$ should be interpreted as a duality between $L^\infty$ and $L^1$.
\end{Rem}

\subsection{Pathwise uniqueness}\label{SS:3.1}
We want to prove that for any pair $(X^{i},\eta^{i},W), i=1,2$, of weak solutions of \eqref{Eq:0.10} defined on the same probability space with the same driving noise $W$ and with $X^{1}_{0} = X^{2}_{0}$, we have $(X^{1},\eta^{1})=(X^{2},\eta^{2})$. This pathwise uniqueness will be used in the next subsection to construct stationary strong solutions of \eqref{Eq:0.10}.
\begin{Prop}\label{Prop:3.1}
Let $x \in \mathcal{C}([0,1],\R^+)$ with $\overline{x} > 0$. Let $(X^{i},\eta^{i},W), i=1,2$ be two weak solutions of \eqref{Eq:0.10} with $X^{1}_{0} = x = X^{2}_{0}$. Then $(X^{1},\eta^{1})=(X^{2},\eta^{2})$.
\end{Prop}
\Proof We use the following Lemma from \cite{MR2349572}. For the sake of completeness, we recall the proof.
\begin{Le}\label{Le:3.1}
Let $\zeta$ be a finite signed measure on $O_{\delta,T}$, $V\in\mathcal{C}(O_{\delta,T})$ and $c> 0$. Suppose that:
\begin{enumerate}
\item[i)] for all $r\in[\delta,T]$, for all $h\in\mathcal{C}([0,1])$, such that $\bar{h}=0$, $\langle h,\zeta\rangle_{O_{r,T}}=0$,
\item[ii)] for all $r\in[\delta,T]$, $\overline{V(r,\cdot)} = c$ with $\langle V,\zeta\rangle_{O_{r,T}}=0$,
\end{enumerate}
then $\zeta$ is the null measure.
\end{Le}
\Proof Let $k\in \mathcal{C}([0,1])$. Since $\zeta$ is a finite measure, by $i)$ we obtain for all $\delta\leq s \leq t \leq T$:
\BD
\langle k,\zeta\rangle_{O_{s,t}} = \langle \bar{k},\zeta\rangle_{O_{s,t}} = \bar{k}\zeta(O_{s,t}), \quad \text{ for all } k \in \mathcal{C}([0,1]).
\ED
This implies $\zeta$ can be decomposed as $\zeta = \gamma \otimes  d\theta$, where $\gamma$ is a measure on $[0,T]$. By $ii)$, we obtain:
\BD
0 = \langle V,\zeta\rangle_{O_{s,t}} = \int_s^t\Big(\int_{0}^1V(s,\theta)d\theta\Big)d\gamma = c\gamma([s,t]).
\ED
We conclude that for all $\delta\leq s \leq t \leq T$, $\gamma([s,t])=0$, since $c>0$. Thus $\zeta$ is the null measure. \EProof
We now prove the proposition. Let $Y(t) = X^1(t,x)-X^2(t,x)$ and $\zeta = \eta^1-\eta^2$, $Y$ is the solution of the following equation:
\BE\label{Eq:3.1}
\left\{\begin{array}{l}
dY=-\frac{1}{2}\left(A^2Y+A(f(X^1)-f(X^2))\right)dt,\\
\\
Y(0)=0.
\end{array}\right.
\EE
We consider now the following approximation of $Y$ :
\BD
Y^N(t,.) = \frac{1}{N}\sum_{n=0}^N\sum_{i=0}^n\langle Y(t),e_i\rangle e_i.
\ED
Since $Y$ is continous, then $Y^N$ converges uniformly to $Y$ on $O_{T}$. Notice that for all $i\geq0$, the process $t \mapsto \langle Y(t),e_{i}\rangle$ has bounded variation, and in particular the process $t \mapsto Y^N(t)$ has bounded variation as process with values in a finite-dimensional subspace of $D(A)$. Taking the scalar product in $V_{-1}$ between $Y$ and $Y^N$, we obtain:
\BE\label{Eq:3.2}
\left\{\begin{array}{l}
d(Y,Y^N)= 2 (Y^N,dY) = -\left(Y^N,A^2Y+A(f(X^1)-f(X^2))+A\zeta\right)dt,\\
\\
(Y,Y^N)(0)=0.
\end{array}\right.
\EE
Moreover for all $t\geq0$, $\langle Y(t), e_{0}\rangle = 0$, so $\Pi Y^N(t) = Y^N(t)$. So we have for all $0< \delta \leq t \leq T$:
\BEA\label{Eq:3.3}
(Y(t),Y^N(t)) &=& (Y(\delta),Y^N(\delta)) +\int_{\delta}^t \langle AY^N(s), Y(s)\rangle ds\nonumber\\
&&+ \int_{\delta}^t \langle Y^N(s),f(X^1(s,x))-f(X^2(s,x))\rangle ds+ \big\langle Y^N(s),\zeta\big\rangle_{O_{\delta,t}}\nonumber\\
&=&(Y(\delta),Y^N(\delta)) - \frac{1}{N}\sum_{n=0}^N\sum_{i=0}^n(i\pi)^2\int_{\delta}^t\langle Y(s),e_i\rangle^2ds+\big\langle Y^N(s),\zeta\big\rangle_{O_{\delta,t}}\nonumber\\
&& +\frac{1}{N}\sum_{n=0}^N\sum_{i=0}^n\int_{\delta}^t\langle Y(s,.),e_i\rangle \langle f(X^1(s,x))-f(X^2(s,x)),e_i\rangle ds\nonumber\\
&&\nonumber\\
&\leq&(Y(\delta),Y^N(\delta)) +\big\langle Y^N(s),\zeta\big\rangle_{O_{\delta,t}}\nonumber\\
&&+\frac{1}{N}\sum_{n=0}^N\sum_{i=0}^n\int_{\delta}^t\langle Y(s,.),e_i\rangle \langle f(X^1(s,x))-f(X^2(s,x)),e_i\rangle ds.
\EEA
For all $s\in[\delta,t]$, 
\BEAS
&&\frac{1}{N}\sum_{n=0}^N\sum_{i=0}^n\langle Y(s),e_i\rangle \langle f(X^1(s,x))-f(X^2(s,x)),e_i\rangle - \langle Y(s), f(X^1(s,x))-f(X^2(s,x))\rangle\\
&=& \langle \frac{1}{N}\sum_{n=0}^N\sum_{i=0}^n\langle Y(s),e_i\rangle e_{i} - Y(s), f(X^1(s,x))-f(X^2(s,x))\rangle\\
&=&\langle Y^N(s) - Y(s), f(X^1(s,x))-f(X^2(s,x))\rangle\\
&\leq& \|Y^N(s) - Y(s)\|_{L^\infty([0,1])}\|f(X^1(s,x))-f(X^2(s,x))\|_{L^1([0,1])},
\EEAS
where $\|\cdot\|_{L^\infty([0,1])}$ and $\|\cdot\|_{L^1([0,1])}$ are the classical norm on the space $[0,1]$. The latter term converges to zero since $Y^N(s)$ converges uniformly to $Y(s)$ on $[0,1]$.
Taking the negative part, we have by Fatou's lemma:
\BEAS
&& \liminf_{N\rightarrow+\infty} \int_{\delta}^t \left(\frac{1}{N}\sum_{n=0}^N\sum_{i=0}^n\langle Y(s),e_i\rangle \langle f(X^1(s,x))-f(X^2(s,x)),e_i\rangle\right)^- ds\\
&\geq&\int_{\delta}^t \liminf_{N\rightarrow+\infty}\left(\frac{1}{N}\sum_{n=0}^N\sum_{i=0}^n\langle Y(s),e_i\rangle \langle f(X^1(s,x))-f(X^2(s,x)),e_i\rangle\right)^- ds\\
&=&\int_{\delta}^t \left(\langle Y(s), f(X^1(s,x))-f(X^2(s,x))\rangle\right)^- ds\\
&=& 0,
\EEAS
since $f$ is nonincreasing.
Taking the limit in \eqref{Eq:3.3} as $N$ grows to infinity, we obtain by the contact condition
\BEAS
\|Y(t)\|_{-1}^2 - \|Y(\delta)\|_{-1}^2&\leq& \big\langle Y,\zeta\big\rangle_{O_{\delta,t}} \\
&= &- \big\langle X^1,\eta^2\big\rangle_{O_{\delta,t}} -\big\langle X^2,\eta^1\big\rangle_{O_{\delta,t}}\\
&\leq&0.
\EEAS
Letting $\delta \rightarrow 0$, we have $Y(t)=0$ for all $t\geq0$ and $X^1(t,x) = X^2(t,x)$ for all $t\geq0$.
Moreover, with the definition of a weak solution, we see that :
\BD
\text{ for all } h\in D(A^2),\qquad  \big\langle Ah,\zeta\big\rangle_{O_{\delta,t}} = 0.
\ED
By density, we obtain $\zeta$ and $V = X^1 = X^2$ satisfy the hypothesis of Lemma \ref{Le:3.1}, and therefore $\zeta$ is the null measure, i.e. $\eta^1=\eta^2$.\EProof

\subsection{Convergence of invariants measures}\label{S:3.2}
Let : 
\BD
K = \{x\in L^2(0,1) ,  x \geq 0\},
\ED
then we know that $\mu_c$ is the law of $Y^c= \B - \overline{\B} +c$. We remark the following inclusion : 
\BD
\{\B_{\theta} \in [-c/2,c/2], \text{ for all } \theta \in [0,1]\} \subset \{Y^c \in K\},
\ED
therefore $\mu_c(K)>0$ with $c>0$. Let us define $U$ the potential associated to the function $f$.
If $f=\fln$ is the logarithmic function, $U$ is defined by:
\BD
U(x)=U_{\ln}(x):=\left\{\begin{array}{rl}
\int_{0}^1\Fln(x(\theta))d\theta &\text{ if } x \in K,\\
+\infty &\text{ else}.
\end{array}\right.
\ED
If $f=f_{\alpha}$ is the negative $\alpha$-power function, $U$ is defined by:
\BD
U(x)=U_{\alpha}(x):=\left\{\begin{array}{rl}
\int_{0}^1F_{\alpha}(x(\theta))d\theta &\text{ if } \int_{0}^1\Big|F_{\alpha}(x(\theta))\Big|d\theta <+\infty \text{ and } x\in K,\\
+\infty &\text{ else}.
\end{array}\right.
\ED
\begin{Rem}\label{Re:3.3}
Note that, for $\alpha<1$, $F_{\alpha}(x(\theta)) = -\frac{1}{1-\alpha}x(\theta)^{1-\alpha}$. By H\"older inequality: 
\BD
\int_{0}^1 |F_{\alpha}(x(\theta))|d\theta <+\infty, \text{ for all } x \in K.
\ED
\end{Rem}
We have the following result : 
\begin{Prop} For $c>0$,
\BD
\nu^n_{c} \rightharpoonup \nu_c := \frac{1}{Z_c}\exp^{-U(x)}\1_{x\in K}\mu_c(dx),\text{ when } n \rightarrow +\infty,
\ED
where $Z_{c}$ is a normalization constant.
\end{Prop}
\Proof
Let $\psi \in \mathcal{C}^0_{b}(L^2,\R)$. We want to prove that 
\BE\label{Eq:3.4}
\int_{H} \psi(x)\exp({-U^n(x)})\mu_{c}(dx) \mathop{\longrightarrow}_{n\rightarrow+\infty} \int_{H} \psi(x)\exp({-U(x)})\1_{x\in K}\mu_{c}(dx).
\EE
\Case[1] $f=\fln$ is the logarithmic function.\\
We have that for a fixed $x\in H$, 
\BE\label{Eq:3.5}
\exp({-U^n(x)}) \mathop{\longrightarrow}_{n\rightarrow+\infty} \exp({-U(x)})\1_{x\in K}.
\EE
Indeed, for all $x\notin K$ there exists $\delta_{x}>0$ small such that $\lambda(\{\theta\in [0,1] / x(\theta)\leq-\delta_{x}\})>0$ and we have:
\BD
\int_{0}^1\Fnln(x(\theta))1_{\{x<0\}}d\theta > \int_{0}^1\Fnln(x(\theta))1_{\{x\leq -\delta_{x}\}}d\theta> 0, \quad\text{for all } n \geq 1.
\ED
Then, since $\Fnln$ is nonincreasing on $(-\infty,0)$:
\BEAS
0\leq\exp({-U^n_{\ln}(x)}) &\leq& \exp\Big(-\int_{0}^1\Fnln(x(\theta))1_{\{x\leq -\delta_{x}\}}d\theta\Big)\\
&\leq& \exp\Big(-\int_{0}^1\Fnln(-\delta_{x})1_{\{x\leq -\delta_{x}\}}d\theta\Big)\\
&\leq& \exp\Big(-\Fnln(-\delta_{x})\lambda(\{x\leq -\delta_{x}\})\Big)\\
&\leq& \exp\Big(\big((1/n -\delta_{x})\ln n - 1 + 1/n\big)\lambda(\{x\leq -\delta_{x}\})\Big).
\EEAS
And this latter term converges to zero as $n$ grows to infinity.\\
Now for $x \in K$, $\Fnln(x(\theta))$ converges to $\Fln(x(\theta))$ almost everywhere as $n$  grows to infinity. Moreover $\Fnln(\tilde{x}(\theta)) \leq \1_{x\leq1} + \Foneln(\tilde{x}(\theta))\1_{x>1}$, and the right-hand side is clearly integrable. By the dominated convergence Theorem, we deduce \eqref{Eq:3.5}.\\
Since $U^n_{\ln} \geq 0$, \eqref{Eq:3.4} follows by dominated convergence Theorem.\\
\noindent\Case[2] $f=f_{\alpha}$ is negative $\alpha$-power function.\\
For a fixed $x \in L^2$, the potentials are increasing as $n$ grows to infinity, we deduce:
\BE\label{Eq:3.6}
\exp({-U^n_{\alpha}(x)}) \leq \exp({-U^1_{\alpha}(x)}), \FAN, \FAX.
\EE
The right-hand side is integrable on $H$, thus it suffices to prove that for a fixed $x\in H$, 
\BE\label{Eq:3.7}
\exp({-U^n_{\alpha}(x)}) \mathop{\longrightarrow}_{n\rightarrow+\infty} \exp({-U_{\alpha}(x)})\1_{x\in K},
\EE
where
\BE
\exp({-U_{\alpha}(x)})\1_{x\in K} = \left\{\begin{array}{ll}
\exp({-U_{\alpha}(x)})& \text{if } \int_{0}^1\left|F_{\alpha}(x(\theta))\right|d\theta <+\infty \text{ and } x\in K,\\
0 & \text{else}.
\end{array}\right.
\EE
For $x\notin K$, there exists $\delta_{x}>0$ small such that $\lambda(\{\theta\in [0,1] / x(\theta)\leq-\delta_{x}\})>0$ and we have:
\BD
\text{for all }n \in \N^*,\qquad \int_{0}^1F^n_{\alpha}(x(\theta))1_{\{x<0\}}d\theta > \int_{0}^1F^n_{\alpha}(x(\theta))1_{\{x\leq -\delta_{x}\}}d\theta> 0.
\ED
Then
\BEAS
0\leq\exp^{-U^n_{\alpha}(x)} &\leq& \exp\left(-\int_{0}^1F^n_{\alpha}(x(\theta))1_{\{x\leq -\delta_{x}\}}d\theta\right)\\
&\leq& \exp\left(-\int_{0}^1F^n_{\alpha}(-\delta_{x})1_{\{x\leq -\delta_{x}\}}d\theta\right)\\
&=& \exp\left(-F^n_{\alpha}\left(-\delta_{x}\right)\lambda\left(\{x\leq -\delta_{x}\}\right)\right).
\EEAS
And this latter term converges to zero as $n$ grows to infinity. Thus \eqref{Eq:3.7} holds.\\
For $x \in K$, such that $\int_{0}^1\left|F_{\alpha}(x(\theta))\right|d\theta <+\infty$, $F^n_{\alpha}(x(\theta))$ converges almost everywhere to $F_{\alpha}(x(\theta))$ as $n$  grows to infinity. Moreover $F^1_{\alpha}(x(\theta)) \leq F^n_{\alpha}(x(\theta)) \leq F_{\alpha}(x(\theta))$ for all $\theta \in  [0,1]$, and by the dominated convergence Theorem \eqref{Eq:3.7} holds.\\
If $\int_{0}^1\left|F_{\alpha}(x(\theta))\right|d\theta =+\infty$, necessarily $\alpha \geq 1$. For $\alpha>1$, $F¬^n_{\alpha} \geq 0$ and \eqref{Eq:3.7} follows from monotone convergence. If $\alpha=1$, we write
\BD
\int_{0}^1 F^n_{\alpha}(x(\theta))d\theta = 
\int_{0}^1 F^n_{\alpha}(x(\theta))\1_{x(\theta)\leq1/2}d\theta + 
\int_{0}^1 F^n_{\alpha}(x(\theta))\1_{x(\theta)>1/2}d\theta. 
\ED
The first term converges to $\int_{0}^1F_{\alpha}(x(\theta))\1_{x(\theta)\leq1/2}$ by monotone convergence, and the second term converges to $\int_{0}^1F_{\alpha}(x(\theta))\1_{x(\theta)<1/2}$ by uniform integrability. We have proved that \eqref{Eq:3.7} always holds, \eqref{Eq:3.4} follows.
\EProof

\subsection{Existence of stationary solutions}\label{SS:3.3}
In this section, we prove the existence of stationary solutions of equation \eqref{Eq:0.10} and that they are limits of stationary solutions of \eqref{Eq:2.1}, in some suitable sense. Fix $c>0$ and consider the unique (in law) stationary solution of \eqref{Eq:2.1} denote $\hat{X}_{c}^{n}$ in $H_{c}$. We are going to prove that the laws of $\hat{X}_{c}^{n}$ weakly converge as $n$ grows to infinity to a stationary strong solution of \eqref{Eq:0.10}.
\begin{Th}\label{Th:3.1}
Let $c>0$ and $T>0$. Almost surely $\hat{X}^n_{c}$ converges as $n$ grows to infinity to a process $\hat{X}_{c}$ in $\mathcal{C}(O_{T})$. Moreover $f(\hat{X}_{c}) \in L^1(O_T)$ almost surely, and setting
\BES
d\eta^n = f^n(\hat{X}^n_{c}(t,\theta))dtd\theta-f(\hat{X}_{c}(t,\theta))dtd\theta,
\EES
then $(\hat{X}^{n}_{c},\eta^n,W)$ converges in law to $(\hat{X}_{c},\eta,W)$ stationary strong solution of \eqref{Eq:0.10}.
\end{Th}
The proof of \ref {Th:3.1} requires arguments that differ significantly in the logarithmic case and in the negative $\alpha$-power case. We thus have chosen to do two separated proofs. Some arguments however are similar and are not repeated.\\
\noindent {\bf{\underline{Proof in the logarithmic case}: }}\\
The proof is splitted in 4 steps. In step 1, assuming that a subsequence of $\hat{X}_{c}^n$ converges in law. Its limit $\hat{X}_{c}$ is shown to satisfy $\fln(\hat{X}_{c}) \in L^1(O_{T})$ almost surely. Then in step 2, under the same assumption as in step 1, we prove that up to a further extraction the measures $\eta^n$ converges to a positive measure $\eta$ and that $(\hat{X}_{c},\eta)$ is a weak solution in the probabilistic sense. It then remains to prove tightness of $\hat{X}_{c}^n$ in step 3 and to use pathwise uniqueness to conclude in step 4.\\
\Step[1] \\
Let us assume that $(n_{k})_{k\in\N}$ is a subsequence such that $(\hat{X}^{n_{k}}_{c})_{n\in \mathbb{N}}$ converges in law in $\mathcal{C}(O_{T})$ to a process $\hat{X}_{c}$.\\
By Skorohod's theorem, we can find a probability space and a sequence of processes $(V^k,\mathcal{W}^{k})_{k\in\N}$ on that probability space such that $(V^k,\mathcal{W}^k)\rightarrow (V,\mathcal{W})$ in $\mathcal{C}(O_{T})$ almost surely and $(V^k,\mathcal{W}^k)$ has the same distribution as $(\hat{X}^{n_{k}}_{c},\mathcal{W})$ for all $k \in \N$. Notice that $V \geq 0$ almost surely since for all $t\leq T$ the law of $V(t,.)$ is $\nu_c$ which is concentrated on $K$. Let now $\xi^k$ and $\rho^k$ be the following measures on $O_{T}$: 
\BD
d\xi^k:=\fln^{n_{k}}(V^k(t,\theta))\1_{V^k< 1}dtd\theta,
\ED
and
\BD
d\rho^k:=\fln^{n_{k}}(V^k(t,\theta))\1_{V^k \geq 1}dtd\theta.
\ED
Let $y \in D(A)$ with $\bar{y}=0$, taking $h \in D(A^2)$ such that $y=Ah$ as a test function in (b) of Definition \ref{De:1}, we deduce that, for all $0\leq t \leq T$,$\big\langle y, \xi^k+\rho^k\big\rangle_{O_{t}}$ has a limit when $n\rightarrow +\infty$.
Moreover by the uniform convergence in $\mathcal{C}(O_{T})$ of $V^k$ to $V$, we have
\BE
\fln^{n_{k}}(V^k(t,\theta))\1_{V^k\geq 1} \mathop{\longrightarrow}_{k\rightarrow+\infty} \fln(V(t,\theta))\1_{V \geq 1},\quad \text{ for all } (t,\theta) \in O_{T},
\EE
and the convergence is uniform. We obtain for all $0\leq t \leq T$ and for all $h \in D(A)$:
\BE\label{Eq:3.10}
\big\langle h, \rho^k\big\rangle_{O_{t}} \mathop{\longrightarrow}_{k\rightarrow+\infty} \int_{O_{t}} h(\theta)\fln(V(s,\theta))\1_{V\geq 1}dsd\theta.
\EE
Note that $\fln(x) \1_{x\geq1}$ is a continuous function so that $\fln(V) \1_{V\geq1} \in L^1(O_{T})$. Moreover, for any $y \in D(A)$ with $\bar{y}=0$, for all $0\leq t\leq T$,
\BE\label{Eq:3.11}
\big\langle y, \xi^k\big\rangle_{O_{t}} \text{ has a limit when } k\rightarrow +\infty.
\EE
Notice that almost surely:
\BE
\fln^{n_{k}}(V^k(t,\theta))\1_{V^k< 1} \mathop{\longrightarrow}_{k\rightarrow+\infty}
\left\{\begin{array}{ll}
+\infty & \text{ if } V(t,\theta)\leq 0,\\
\fln(V(t,\theta)) & \text{ if } V(t,\theta) \in (0,1].
\end{array}\right.
\EE
Thus the limit of this term is not trivial. Let us now prove that the total mass $\xi^n(O_{T})$ is bounded. We use the following Lemma whose proofs is postponed to the end of this section.
\begin{Le}\label{Le:3.2}
Let $T>0$, and $\{\mu^k\}_{k\in\N}$ be a sequence of finite positive measures on $O_{T}$. Suppose there exists $\{w^k\}_{k\in\N}$ a sequence of functions in $\mathcal{C}(O_{T})$ such that $w^k$ converges uniformly to $w$, when $k$ grows to infinity. Suppose also there exist a function $M_{T} : \mathcal{C}(O_{T}) \rightarrow \R^+$ and two nonnegative constants  $m_{T}$ and $c_{T}$ such that 
\BE\label{Eq:3.12}
\text{for all } h\in D(A) \text{ such that } \bar{h}=0,\quad\big\langle h, \mu^k\big\rangle_{O_{T}}\leq M_{T}(h),\quad \text{ for all } k\in\N,
\EE
\BE
\text{for all } t \in O_{T}, \int_{0}^{1}w(t,\theta)d\theta = c_{T} >0
\EE
and
\BE\label{Eq:3.14}
\big\langle w^k, \mu^k\big\rangle_{O_{T}} \leq m_{T}.
\EE
Then there exists a constant $\tilde{M}_{T}$ such that
\BE\label{Eq:3.15}
\text{ for all } h \in \mathcal{C}(O_{T}),\quad\big\langle h, \mu^k\big\rangle_{O_{T}} \leq \tilde{M}_{T}\| h\|_{\infty},\quad  \text{for all } k\in\N.
\EE
and in particular $\mu^k(O_{T})$ is bounded uniformly for $k\in\N$.
\end{Le}
Let us denote by :
\BE
M_{T}(h) = \sup_{k\in \N}\left|\big\langle h, \xi^k\big\rangle_{O_{T}}\right|
\EE
for $h \in D(A)$ such that $\bar{h}=0$. By \eqref{Eq:3.11}, we know that $M_{T}$ is well defined.
Moreover we have
\BE
\big\langle (V^k)^+, \xi^k\big\rangle_{O_{T}} = \int_{O_{T}} (V^k(t,\theta))^+\fln^{n_{k}}(V^k(t,\theta))\1_{V^k<1}dtd\theta.
\EE
Since $(x)^+\fln^{n_{k}}(x)\1_{x<1}$ is uniformly bounded in $k\in\N$, there exists a positive constant $m_{T}$ such that 
\BE
\big\langle (V^k)^+, \xi^k\big\rangle_{O_{T}} \leq m_{T}.
\EE
Since $V$ is almost surely positive, $(V^k)^+$ converges uniformly to $V$ and $\overline{V(t,.)}=c_{T}>0$ for all $t \in [0,T]$. We use Lemma \ref{Le:3.2} and obtain $\limsup_{k\rightarrow+\infty}\xi^k(O_{T}) < +\infty$.\\
Thanks to Fatou Lemma, we can write : 
\BEA\label{Eq:3.27}
 \int_{O_{T}} \big[\fln(V(s,\theta))\1_{V<1} \big]ds d\theta &= &\int_{O_{T}}\mathop{\liminf}_{k\rightarrow+\infty} \big[\fln^{n_{k}}(V^k(s,\theta))\1_{V^n<1}ds d\theta \big]\nonumber\\
&\leq&\mathop{\liminf}_{k\rightarrow+\infty} \int_{O_{T}}\big[\fln^{n_{k}}(V^k(s,\theta))\1_{V^k<1}\big]ds d\theta \\
& <& +\infty.\nonumber
\EEA
It follows that almost surely $\fln(V) \in L^1(O_{T})$.\EStep[1]
\indent\Step[2] \\
We again assume that we have $(n_{k})_{k\in\N}$ a subsequence such that $(\hat{X}^{n_{k}}_{c})_{k\in \mathbb{N}}$ converges in law to a process $\hat{X}_{c}$. Again, by Skorohod's theorem, we can find a probability space and a sequence of processes $(V^k,\mathcal{W}^k)_{k\in\N}$ such that almost surely $(V^k,\mathcal{W}^k)\rightarrow (V,\mathcal{W})$ in $\mathcal{C}(O_{T})$ as $k$ grows to infinity, and $(V^k,\mathcal{W}^k)$ has the same distribution as $(\hat{X}^{n_{k}}_{c},W)$ for all $k \in \N$.\\
By step 1, the total mass $\xi^k(O_{T})$ is bounded and there exists $(n_{k_{m}})_{m\in\N}$ a sub-subsequence such that the measures
\BD
\xi^{k_{m}}:=\fln^{n_{k_{m}}}(V^{k_{m}}(t,\theta))\1_{V^{k_{m}}< 1}dtd\theta
\ED
converge to a measure $\xi$. \\
We denote by $\lambda$ the following measure:
\BE
d\lambda := \fln(V(t,\theta))\1_{V<1}dtd\theta,
\EE
and $\zeta^m := \xi^{k_{m}} -\lambda$. Thus $\zeta^m$ converges to the measure $\zeta:=\xi -\lambda$. Let $u$ be a continuous nonnegative function on $O_{T}$, we have
\BEAS
\big\langle u, \zeta\big\rangle_{O_{T}}& =& \lim_{m\rightarrow+\infty} \big\langle u, \zeta^m\big\rangle_{O_{T}}\\
& = &\lim_{m\rightarrow+\infty} \big\langle u, \xi^{k_{m}}\big\rangle_{O_{T}} -\big\langle u, \lambda\big\rangle_{O_{T}}.
\EEAS
And this is positive, thanks to \eqref{Eq:3.27}. Therefore $\zeta$ is a positive measure. Taking the limit as $m$ grows to infinity in the approximated equation, we obtain that for all $h \in D(A^2)$ and for all $0\leq t \leq T$:
\BEAS
\langle V(t,.),h\rangle&=&\langle x,h\rangle - \int_{O_{t}} V(s,\theta)A^2h(\theta) dsd\theta - \int_{O_{t}} \fln(V(s,\theta))Ah(\theta) dsd\theta\\
&& -\big\langle Ah, \zeta\big\rangle_{O_{t}} -  \int_{0}^{t}\langle Bh, d\mathcal{W}\rangle.
\EEAS
This is the expected equation. Let us now show that the contact condition holds for $(V,\zeta)$. We prove in fact that for all $\beta$ non negative:
\BE\label{Eq:3.28}
0\leq\big\langle V, \zeta\big\rangle_{O_{T}}\leq\beta.
\EE
The key is to study the behavior of $\fln^{n_{k_{m}}}(V^{k_{m}}(t,\theta))\1_{V^{k_{m}}<1}$ near points $(t,\theta)\in O_{T}$ such that $V(t,\theta)$ is small. Fix $\beta>0$, there exists $\varepsilon>0$ such that $-T\varepsilon\ln(\varepsilon)\leq\beta$. Let us define the following measures for all $m\in\N$.
\BES
\begin{array}{ll}
d\xi^{m}_{\varepsilon} := \fln^{n_{k_{m}}}(V^{k_{m}}(t,\theta))\1_{V^{k_{m}}<\varepsilon}dtd\theta,
&d\tau^{m}_{\varepsilon} := \fln^{n_{k_{m}}}(V^{k_{m}}(t,\theta))\1_{\varepsilon\leq V^{k_{m}}<1}dtd\theta,\\
\\	
d\lambda_{\varepsilon} := \fln(V(t,\theta))\1_{V<\varepsilon}dtd\theta, &d\tau_{\varepsilon} := \fln(V(t,\theta))\1_{\varepsilon\leq V<1}dtd\theta.
\end{array}
\EES
Clearly $\tau_{\varepsilon}^m$ converges to $\tau_{\varepsilon}$, it follows
\BEAS
\limsup_{m\rightarrow +\infty}\  \langle V^{k_{m}}, \zeta^m\big\rangle_{O_{T}}&=&\limsup_{m\rightarrow +\infty}\left( \big\langle V^{k_{m}},\xi^{k_{m}}_{\varepsilon}\big\rangle_{O_{T}}\!\!-\big\langle V^{k_{m}},\lambda_{\varepsilon}\big\rangle_{O_{T}} +  \big\langle V^{k_{m}},\tau^m_{\varepsilon}\big\rangle_{O_{T}}\!\!-\big\langle V^{k_{m}},\tau_{\varepsilon}\big\rangle_{O_{T}}\right)\\
\\
&=&\limsup_{m\rightarrow +\infty}\left(\int_{O_{T}} \!\!V^{k_{m}}\fln^{n_{k_{m}}}(V^{k_{m}})\1_{V^{k_{m}}<\varepsilon}dtd\theta - \int_{O_{T}} \!\!V^{k_{m}}\fln(V)\1_{V<\varepsilon}dtd\theta\right)\\
\\
&\leq&\limsup_{m\rightarrow +\infty}	\left(\int_{O_{T}} V^{k_{m}}\fln^{n_{k_{m}}}(V^{k_{m}})\1_{0\leq V^{k_{m}}<\varepsilon}dtd\theta\right)\\
&& + \limsup_{m\rightarrow +\infty}	\left(\int_{O_{T}} (V^{k_{m}})^-\fln(V)\1_{V<\varepsilon}dtd\theta \right)
\EEAS
Since $(V^{k_{m}})^-$ converges uniformly to zero, we deduce:
\BEAS
\limsup_{m\rightarrow +\infty} \ \langle V^{k_{m}}, \zeta^m\big\rangle_{O_{T}}&\leq&T\limsup_{m\rightarrow +\infty}\sup_{x\in[0,\varepsilon]}\left|-x\ln\left(x+\frac{1}{n_{k_{m}}}\right)\right|\\
\\
&\leq&T\limsup_{m\rightarrow +\infty}\left(-\varepsilon\ln\left(\varepsilon+\frac{1}{n_{k_{m}}}\right)\right)\\
\\
&\leq&-T\varepsilon\ln\left(\varepsilon\right).
\EEAS
Thus the contact condition holds.\EStep[2]
\indent\Step[3]\\
By the convergence of the family $(\nu^n_{c})_{n\in \N}$, we know that the initial distribution of $\hat{X}_{c}^{n}$ converges to $\nu_c$. We now follow the same argument as in \cite{MR2349572}, to prove for all $T>0$, the laws of $(\hat{X}^{n}_{c})_{n\in \N}$ are tight in $\mathcal{C}(O_{T})$. Fix $n\geq 1$ and $T>0$, by the Lyons-Zheng's decomposition (see theorem 5.7.1 in \cite{MR1303354}), we can find $M^h$, respectively $N^h$, two orthogonal martingales with respect to the natural filtration of $\left(\hat{X}_{c}^n(t)\right)_{0\leq t\leq T}$, respectively the natural filtration of $\left(\hat{X}_{c}^n(T-t)\right)_{0\leq t\leq T}$, such that for all $t\in [0,T]$ and for all $h\in H$:
\BE\label{Eq:3.30a}
(h,\hat{X}^{n}_{c}(t)-\hat{X}^{n}_{c}(0)) = \frac{1}{2} M^h_{t} - \frac{1}{2}(N^h_{T}-N^h_{T-t}).
\EE
Moreover, the quadratic variations are both equal to $\langle\langle M^h\rangle\rangle_{t} = \langle\langle N^h\rangle\rangle_{t} = t \|\Pi h\|_{-1}^2$. Let $u_{i} = Q^{\gamma-1} e_{i}$ for all $i\geq 1$. To simplify the notations, we denote $M^i$ and $N^i$ the martingales $M^{u_{i}}$ and $N^{u_{i}}$ defined in \eqref{Eq:3.30a}. Then we have:
\BEAS
\E\left[\left\|\hat{X}_c^{n}(t)-\hat{X}_c^{n}(s)\right\|^2_{-\gamma}\right]
&=&\E\left[\sum_{i=1}^{+\infty} \left((\hat{X}_c^{n}(t)-\hat{X}_c^{n}(s),v_{i})\right)^2_{-\gamma}\right]\\
&=&\E\left[\sum_{i=1}^{+\infty} \left((\hat{X}_c^{n}(t)-\hat{X}_c^{n}(0),u_{i})_{-1}- (\hat{X}_c^{n}(s)-\hat{X}_c^{n}(0),u_{i})\right)^2_{-1}\right]\\
&=&\frac{1}{4}\E\left[\sum_{i=1}^{+\infty} \left(M^i_{t}+N^i_{T-t}-M^i_{s}-N^i_{T-s}\right)^2_{-1}\right]\\
&=&\frac{1}{4}\sum_{i=1}^{+\infty}\E\left[\langle\langle M^i\rangle\rangle_{t-s}+\langle\langle N^i\rangle\rangle_{t-s} \right]\\
&=&\frac{1}{2}\sum_{i=1}^{+\infty}|t-s|\|u_{i}\|^2_{-1}\\
&=&\frac{|t-s|}{2}\sum_{i=1}^{+\infty}(i\pi)^{2-2\gamma}\\
&\leq&\frac{|t-s|}{2} K^2_{-\gamma},
\EEAS
where $K_{-\gamma}$ is the Hilbert-Schmidt's norm of inclusion of $V_{-1}$ in $V_{-\gamma}$ which is finite for $\gamma>3/2$.
So we have found  a constant $C>0$ such that for all $t,s\in [0,T]$:
\BE\label{Eq:3.31a}
\left(\E\left[\left\|\hat{X}_c^{n}(t)-\hat{X}_c^{n}(s)\right\|^2_{-\gamma}\right]\right)^{\frac{1}{2}}\leq C|t-s|^{\frac{1}{2}}.
\EE
Furthermore, for $0<\delta<1/2$ and $r\geq 1$, since $\hat{X}_{c}^n$ is a stationary solution, there exists $C$ such that for all $t,s\in [0,T]$:
\BEA\label{Eq:3.32a}
\left(\E\left[\left\|\hat{X}_c^{n}(t)-\hat{X}_c^{n}(s)\right\|^2_{W^{\delta,r}(0,1)}\right]\right)^{\frac{1}{2}}
&\leq& \left(\E\left[\left\|\hat{X}_c^{n}(t)\right\|^2_{W^{\delta,r}(0,1)}\right]\right)^{\frac{1}{2}}+ \left(\E\left[\left\|\hat{X}_c^{n}(s)\right\|^2_{W^{\delta,r}(0,1)}\right]\right)^{\frac{1}{2}}\nonumber\\
&\leq& 2\left(\int_H\|y\|^2_{W^{\delta,r}(0,1)}\nu_c^n(dy)\right)^{\frac{1}{2}}\nonumber\\
&\leq& C\left(\int_H\|y\|^2_{W^{\delta,r}(0,1)}\mu_c(dy)\right)^{\frac{1}{2}},
\EEA
since $U^n \geq 0$. And this latter term is finite.\\
Let $\kappa \in [0,1]$ and set $\lambda=\kappa \delta - (1-\kappa)\gamma$, $\frac{1}{q}=\kappa \frac{1}{r}+(1-\kappa)\frac{1}{2}$. Then by interpolation for all $t,s\in [0,T]$:
\BEAS
\Big(\E\Big[\big\|\hat{X}_c^{n}(t)&-&\hat{X}_c^{n}(s)\big\|^2_{W^{\lambda,q}(0,1)}\Big]\Big)^{\frac{1}{2}}\\
&\leq&\left(\E\left[\left\|\hat{X}_c^{n}(t)-\hat{X}_c^{n}(s)\right\|^2_{W^{\delta,r}(0,1)}\right]\right)^{\frac{\kappa}{2}}
\left(\E\left[\left\|\hat{X}_c^{n}(t)-\hat{X}_c^{n}(s)\right\|^2_{-\gamma}\right]\right)^{\frac{(1-\kappa)}{2}}.
\EEAS
We use \eqref{Eq:3.31a}, \eqref{Eq:3.32a} and Sobolev embedding to conclude. Indeed, for any $\beta \in (0,1/2)$, we can choose $\delta \in (0,1/2)$, $\gamma >2$, $r\geq1$ and $\kappa \in (0,1)$ such that $(\lambda-\beta)q>1$. It follows that there exists $C$ such that for all $s,t \in [0,T]$:
\BD
\left(\E\left[\left\|\hat{X}_c^{n}(t)-\hat{X}_c^{n}(s)\right\|^2_{\mathcal{C}^{\beta}([0,1])}\right]\right)^{\frac{1}{2}}\leq C |t-s|^{\frac{1-\kappa}{4}},
\ED
Finally, we can conclude by the theorem 7.2 in chapter 3 of \cite{MR838085} that the laws of $(\hat{X}^{n}_{c})_{n\in \N}$ are tight in $\mathcal{C}(O_{T})$.\EStep[3]
\indent\Step[4]\\
We use a Lemma in \cite{MR1392450}. This lemma allows to get the convergence of the approximated solutions in probability in any space in which these approximated solutions are tight.
\begin{Le}\label{Le:3.3}
Let $\{Z_{n}\}_{n\geq1}$ be a sequence of random elements on a Polish space $E$ endowed by its borel $\sigma$-algebra. Then $\{Z_{n}\}_{n\geq1}$ converges in probability to an $E$-valued 
random element if and any if from every pair of subsequences $\{(Z_{n_{k}^1},Z_{n_{k}^2})_{k\geq1}$, one 
can extract a subsequence which converges weakly to a random element supported on the diagonal 
$\{(x, y) \in E \times E, x = y\}$.
\end{Le}
For any subsequence $(n_{k})_{k\in\N}$, we have convergence of $\xi^{k}$ to a finite measure $\xi$ on $O_{T}$ along some sub-subsequence $(k_{m})_{m\in\N}$. Let $\xi_i, i=1..2$ be two such limits. By the second step, and the uniqueness of the reflexion measure, we know $\zeta_{1} := \xi_{1}-\lambda$ and $\zeta_{2} := \xi_{2}-\lambda$ are equals. So the limit of $(\xi^k)_{k\in\N}$ is unique, and $\xi^k$ converges to its limit $\xi$.\\
Assume $(n^1_{k})_{k\in\N}$ and $(n^1_{k})_{k\in\N}$ are two arbitrary subsequences. In the notations of the second step and by the third step, the process $\left(\hat{X}^{n^1_{k}}_{c},\hat{X}^{n^2_{k}}_{c},W\right)$ is tight in a suitable space. By Skorohod's therorem, we can find a probability space and a sequence of processes $(V_1^k,V_2^k,\mathcal{W}^k)$ such that $(V_1^k,V_2^k,\mathcal{W}^k)\rightarrow(V_{1},V_{2},\mathcal{W})$ almost surely in $\mathcal{C}(O_{T})$, and $(V_1^k,V_2^k,\mathcal{W}^k)$ as the same distribution as $\left(\hat{X}^{n^1_{k}}_{c},\hat{X}^{n^2_{k}}_{c},W\right)$ for all $k\in\N$. In the Skorohod's space, the approximated measures respectively converge to two contact measures $\zeta_{1}$ and $\zeta_{2}$. By the second step, $(V_{1},\zeta_{1},\mathcal{W})$ and $(V_{2},\zeta_{2},\mathcal{W})$ are both weak solutions of \eqref{Eq:0.10}. By uniqueness, necessarily $V_{1}=V_{2}$ and $\zeta_{1}=\zeta_{2}$. Therefore the subsequence $\left(\left(\hat{X}^{n^1_{k}}_{c},\eta^{n_{k}^{1}},W\right),\left(\hat{X}^{n^2_{k}}_{c},\eta^{n_{k}^{2}},W\right)\right)_{k\in\N}$ converges in law to a process supported on the diagonal. We use Lemma \ref{Le:3.3} to prove that the sequence $(\hat{X}^{n}_{c},\eta^n,W)$ converges in law to $(\hat{X}_{c},\eta,W)$ stationary weak solution of \eqref{Eq:0.10}. Moreover by pathwise uniqueness and existence of strong solutions, we obtain that every weak solution is also a strong solution.\EStep[4]
\noindent {\bf{\underline{Proof in negative $\alpha$-power case}: }}\\
We again split the proof in four steps.\\
\Step[1]\\
Let us assume that $(n_{k})_{k\in\N}$ is a subsequence such that $(\hat{X}^{n_{k}}_{c})_{n\in \mathbb{N}}$ converges in law in $\mathcal{C}(O_{T})$ to a process $\hat{X}_{c}$.\\
By Skorohod's theorem, we can find a probability space and a sequence of processes $(V^k,\mathcal{W}^{k})_{k\in\N}$ on that probability space such that $(V^k,\mathcal{W}^k)\rightarrow (V,\mathcal{W})$ in $\mathcal{C}(O_{T})$ almost surely and $(V^k,\mathcal{W}^k)$ has the same distribution as $(\hat{X}^{n_{k}}_{c},\mathcal{W})$ for all $k \in \N$. Notice that $V \geq 0$ almost surely since for all $t\leq T$ the law of $V(t,.)$ is $\nu_c$ which is concentrated on $K$. Let now $\xi^k$ be the following measure on $O_{T}$: 
\BD
d\xi^k:=f^{n_{k}}_{\alpha}(V^k(t,\theta))dtd\theta.
\ED
Let $y \in D(A)$ with $\bar{y}=0$, taking $h \in D(A^2)$ such that $y=Ah$ as a test function in (b) of Definition \ref{De:1}, we deduce that, for all $0\leq t \leq T$,
\BE\label{Eq:3.30}
\big\langle y, \xi^k\big\rangle_{O_{t}} \text{ has a limit when } k\rightarrow +\infty.
\EE
Like in the logarithmic case, we now prove that the total mass $\xi^n(O_{T})$ is bounded. Let us denote by :
\BE
M_{T}(h) = \sup_{k\in \N}\left|\big\langle h, \xi^k\big\rangle_{O_{T}}\right|
\EE
for $h \in D(A)$ such that $\bar{h}=0$. By \eqref{Eq:3.30}, we know that $M_{T}$ is well defined.
Therefore by Lemma \ref{Le:3.2} it suffices to find a function $w$, such that $\overline{w(t,.)} = c_{T} >0$, and a sequence $(w^k)_{k\in\N}$ such that for a positive constant $m_{T}$, $w^k$ converges uniformly to $w$ such that
\BE\label{Eq:3.31}
\big\langle w^k, \xi^k\big\rangle_{O_{T}} \leq m_{T}.
\EE
Denote by $w^k := ((V^k)^++1/n_{k})^{\alpha}$, it converges uniformly to $w:=V^\alpha$. Since $w^k f^{n_{k}}_{\alpha}(V^k(t,\theta))=1$, \eqref{Eq:3.31} holds with $m_{T}=T$. As in the logarithmic case, by Fatou Lemma, it follows that almost surely $f_{\alpha}(V) \in L^1(O_{T})$.\EStep[1]
\indent\Step[2]
We again assume that we have $(n_{k})_{k\in\N}$ a subsequence such that $(\hat{X}^{n_{k}}_{c})_{k\in \mathbb{N}}$ converges in law to a process $\hat{X}_{c}$. Again, by Skorohod's theorem, we can find a probability space and a sequence of processes $(V^k,\mathcal{W}^k)_{k\in\N}$ such that almost surely $(V^k,\mathcal{W}^k)\rightarrow (V,\mathcal{W})$ in $\mathcal{C}(O_{T})$ as $k$ grows to infinity, and $(V^k,\mathcal{W}^k)$ has the same distribution as $(\hat{X}^{n_{k}}_{c},W)$ for all $k \in \N$.\\
By step 1, the total mass $\xi^k(O_{T})$ is bounded and there exists $(n_{k_{m}})_{m\in\N}$ a sub-subsequence such that the measures
\BD
\xi^{k_{m}}:=f^{n_{k_{m}}}_{\alpha}(V^{k_{m}}(t,\theta))\1_{V^{k_{m}}< 1}dtd\theta
\ED
converges to a measure $\xi$. \\
We denote by $\lambda$ the following measure:
\BE
d\lambda := f_{\alpha}(V(t,\theta))dtd\theta,
\EE
and $\zeta^m := \xi^m -\lambda$. Thus $\zeta^m$ converges to the measure $\zeta:=\xi -\lambda$. Thanks to Fatou Lemma, $\zeta$ is a positive measure. Taking the limit as $m$ grows to infinity in the approximated equation, we obtain that for all $h \in D(A^2)$ and for all $0\leq t \leq T$:
\BEAS
\langle V(t,.),h\rangle&=&\langle x,h\rangle - \int_{O_{t}} V(s,\theta)A^2h(\theta) dsd\theta - \int_{O_{t}} f_{\alpha}(V(s,\theta))Ah(\theta) dsd\theta \\
&&-\big\langle Ah, \zeta\big\rangle_{O_{t}} -  \int_{0}^{t}\langle Bh, d\mathcal{W}\rangle.
\EEAS
This is the expected equation. Let us now show that the contact condition holds for $(V,\zeta)$.\\
\Case[1] : $0\leq\alpha<1$.\\
As in the second step of the logarithmic case, fix $\beta>0$, so there exists $\varepsilon>0$ such that $T\varepsilon^{1-\alpha}\leq\beta$. Let us define the following measures for all $m\in\N$:
\BES
\begin{array}{ll}
d\xi^{m}_{\varepsilon} := f_{\alpha}^{n_{k_{m}}}(V^{k_{m}}(t,\theta))\1_{V^{k_{m}}<\varepsilon}dtd\theta,&
d\tau^m_{\varepsilon} := f_{\alpha}^{n_{k_{m}}}(V^{k_{m}}(t,\theta))\1_{\varepsilon\leq V^{k_{m}}}dtd\theta,\\
\\
d\lambda_{\varepsilon} := f_{\alpha}(V(t,\theta))\1_{V<\varepsilon}dtd\theta,& d\tau_{\varepsilon} := f_{\alpha}(V(t,\theta))\1_{\varepsilon\leq V}dtd\theta.
\end{array}
\EES
Since $\tau_{\varepsilon}^m$ converges to $\tau_{\varepsilon}$, we have
\BEAS
\limsup_{m\rightarrow +\infty} \ \langle V^{k_{m}}, \zeta^m\big\rangle_{O_{T}} &=&\limsup_{m\rightarrow +\infty}\left( \big\langle V^{k_{m}},\xi^{k_{m}}_{\varepsilon}\big\rangle_{O_{T}}\!\!-\big\langle V^{k_{m}},\lambda_{\varepsilon}\big\rangle_{O_{T}} +  \big\langle V^{k_{m}},\tau^m_{\varepsilon}\big\rangle_{O_{T}}\!\!-\big\langle V^{k_{m}},\tau_{\varepsilon}\big\rangle_{O_{T}}\right)\\
\\
&=&\limsup_{m\rightarrow +\infty}\left(\int_{O_{T}} \!\!V^{k_{m}}f_{\alpha}^{n_{k_{m}}}(V^{k_{m}})\1_{V^{k_{m}}<\varepsilon}dtd\theta - \int_{O_{T}} \!\!V^{k_{m}}f_{\alpha}(V)\1_{V<\varepsilon}dtd\theta\right)\\
\\
&\leq&\limsup_{m\rightarrow +\infty}	\left(\int_{O_{T}} V^{k_{m}}f_{\alpha}^{n_{k_{m}}}(V^{k_{m}})\1_{0\leq V^{k_{m}}<\varepsilon}dtd\theta\right)\\
&& + \limsup_{m\rightarrow +\infty}	\left(\int_{O_{T}} (V^{k_{m}})^-f_{\alpha}(V)\1_{V<\varepsilon}dtd\theta \right).
\EEAS
It follows
\BEAS
\limsup_{m\rightarrow +\infty} \ \langle V^{k_{m}}, \zeta^m\big\rangle_{O_{T}}&\leq&T\limsup_{m\rightarrow +\infty}\sup_{x\in[0,\varepsilon]}\left|\frac{x}{\left(x+\frac{1}{n_{k_{m}}}\right)^\alpha}\right|\\
\\
&\leq&T\limsup_{m\rightarrow +\infty}\left|\frac{\varepsilon}{\left(\varepsilon+\frac{1}{n_{k_{m}}}\right)^\alpha}\right|\\
\\
&\leq&T \varepsilon^{1-\alpha}.
\EEAS
Thus the contact condition holds.\\
\Case[2] : $\alpha\geq1$.\\
Let $\gamma>0$, we prove that for all nonegative $\beta$, $0\leq\big\langle V^{\alpha+\gamma}, \zeta\big\rangle_{O_{T}}\leq\beta$ and conclude that the contact condition holds by H\"older inequality.\\
Fix $\beta>0$, so there exists $\varepsilon>0$ such that $T\varepsilon^{\gamma}\leq\beta$. Let us define the following measures for all $m\in\N$:
\BES
\begin{array}{ll}
d\xi^{m}_{\varepsilon} := f_{\alpha}^{n_{k_{m}}}(V^{k_{m}}(t,\theta))\1_{V^{k_{m}}<\varepsilon}dtd\theta,&
d\tau^m_{\varepsilon} := f_{\alpha}^{n_{k_{m}}}(V^{k_{m}}(t,\theta))\1_{\varepsilon\leq V^{k_{m}}}dtd\theta,\\
\\
d\lambda_{\varepsilon} := f_{\alpha}(V(t,\theta))\1_{V<\varepsilon}dtd\theta, &
d\tau_{\varepsilon} := f_{\alpha}(V(t,\theta))\1_{\varepsilon\leq V}dtd\theta.
\end{array}
\EES
Since $\tau_{\varepsilon}^m$ converges to $\tau_{\varepsilon}$, we have
\BEAS
\limsup_{m\rightarrow +\infty} \ \langle \left(V^{k_{m}}\right)^{\alpha+\gamma}, \zeta^m\big\rangle_{O_{T}} &=&\limsup_{m\rightarrow +\infty}\Big( \big\langle \left(V^{k_{m}}\right)^{\alpha+\gamma},\xi^{k_{m}}_{\varepsilon}\big\rangle_{O_{T}}-\big\langle \left(V^{k_{m}}\right)^{\alpha+\gamma},\lambda_{\varepsilon}\big\rangle_{O_{T}} \\
\\
&&+  \big\langle \left(V^{k_{m}}\right)^{\alpha+\gamma},\tau^m_{\varepsilon}\big\rangle_{O_{T}}-\big\langle \left(V^{k_{m}}\right)^{\alpha+\gamma},\tau_{\varepsilon}\big\rangle_{O_{T}}\Big)\\
\\
&=&\limsup_{m\rightarrow +\infty}\Bigg(\int_{O_{T}} \left(V^{k_{m}}\right)^{\alpha+\gamma}f_{\alpha}^{n_{k_{m}}}(V^{k_{m}})\1_{V^{k_{m}}<\varepsilon}dtd\theta \\
\\
&&- \int_{O_{T}} \left(V^{k_{m}}\right)^{\alpha+\gamma}f_{\alpha}(V)\1_{V<\varepsilon}dtd\theta\Bigg)\\
\\
&\leq&\limsup_{m\rightarrow +\infty}	\left(\int_{O_{T}} \left(V^{k_{m}}\right)^{\alpha+\gamma}f_{\alpha}^{n_{k_{m}}}(V^{k_{m}})\1_{0\leq V^{k_{m}}<\varepsilon}dtd\theta\right)\\
&& + \limsup_{m\rightarrow +\infty}	\left(\int_{O_{T}} \left(\left(V^{k_{m}}\right)^-\right)^{\alpha+\gamma}f_{\alpha}(V)\1_{V<\varepsilon}dtd\theta \right).
\EEAS
It follows
\BEAS
\limsup_{m\rightarrow +\infty} \ \langle \left(V^{k_{m}}\right)^{\alpha+\gamma}, \zeta^m\big\rangle_{O_{T}}
&\leq&T\limsup_{m\rightarrow +\infty}\sup_{x\in[0,\varepsilon]}\left|\frac{x^{\alpha + \gamma}}{\left(x+\frac{1}{n_{k_{m}}}\right)^\alpha}\right|\\
\\
&\leq&T\limsup_{m\rightarrow +\infty}\left|\frac{\varepsilon^{\alpha + \gamma}}{\left(\varepsilon+\frac{1}{n_{k_{m}}}\right)^\alpha}\right|\\
\\
&\leq&T \varepsilon^{\gamma}.
\EEAS
Thus the contact condition holds.
\EStep[2]
Step 3 and step 4 are strictly identical to the logarithmic case and we do not repeat them. This ends the proof of Theorem \ref{Th:3.1}. Now we give the proof of the Lemma \ref{Le:3.2}.\\
\noindent {\bf{\underline{Proof of the Lemma \ref{Le:3.2}}: }}\\
We prove this Lemma thanks to the previous Lemma \ref{Le:3.1}. If $\mu^k(O_{T})$ is bounded uniformly for $k\in\N$, then the constant
\BE\label{Eq:3.16}
\tilde{M}_{T}=\sup_{k\in\N}\mu^k(O_{T}) 
\EE
satisfies \eqref{Eq:3.15}. Suppose $\mu^k(O_{T})$ is unbounded, then there exists $k_{0} \in\N$ such that $\mu^k(O_{T})>0$ for all $k\geq k_{0}$, we denote for all $k \geq k_{0}$
\BES
\nu^k := \mu^k/\mu^k(O_{T}).
\EES
$\{\nu^k\}_{k\geq k_{0}}$ is a sequence of probability measure on $O_{T}$, and we can extract a subsequence $\{\nu^{k_{m}}\}_{m\in\N}$ such that there exists a probability measure $\nu$ with $\nu^{k_{m}} \rightharpoonup \nu$ when $m$ grows to infinity. Therefore, by the uniform convergence of $w^k$
\BE
\big\langle w^{k_{m}}, \nu^{k_{m}}\big\rangle_{O_{T}}  \mathop{\longrightarrow}_{l\rightarrow +\infty} \big\langle w, \nu\big\rangle_{O_{T}}.
\EE
And by the uniform boundedness in \eqref{Eq:3.14}, we have
\BE
\big\langle w^{k_{m}}, \nu^{k_{m}}\big\rangle_{O_{T}} \leq \frac{m_{T}}{\mu^{k_{m}}(O_{T})} \mathop{\longrightarrow}_{l\rightarrow +\infty} 0,
\EE
therefore
\BE
\big\langle w, \nu\big\rangle_{O_{T}} = 0.
\EE
Moreover, for all $h\in D(A)$ such that $\bar{h}=0$
\BE
\big\langle h, \nu^{k_{m}}\big\rangle_{O_{T}} \mathop{\longrightarrow}_{l\rightarrow +\infty} \big\langle h, \nu\big\rangle_{O_{T}},
\EE
and by the uniform boundedness in \eqref{Eq:3.12}, for all $h\in D(A)$ such that $\bar{h}=0$, we have
\BE
\big\langle h, \nu^{k_{m}}\big\rangle_{O_{T}} \leq \frac{M_{T}(h)}{\mu^{k_{m}}(O_{T})} \mathop{\longrightarrow}_{l\rightarrow +\infty} 0.
\EE
So that for all $h\in D(A)$ such that $\bar{h}=0$, we have
\BE\label{Eq:3.37}
\big\langle h, \nu\big\rangle_{O_{T}} = 0.
\EE
Since $\nu$ is a probability measure, we deduce that \eqref{Eq:3.37} holds in fact for any $h \in\mathcal{C}(O_{T})$ such that $\bar{h}=0$. The hypothesis of Lemma \ref{Le:3.1} are satisfied, and we can conclude that the measure $\nu$ is null. This is a contradiction since $\nu$ is a probability measure. Then the sequence $\mu^k(O_{T})$ is bounded uniformly for $k\in\N$, and the constant $\tilde{M}_{T}$ in \eqref{Eq:3.16} fulfills \eqref{Eq:3.15}.
\EProof

\subsection{Convergence of the semigroup}\label{SS:3.4}
First we state the following result which is a corollary of Theorem \ref{Th:3.1}.
\begin{Coro}\label{Co:3.1}Let $c>0$.
\begin{enumerate}
\item[i)] There exists a continuous process $(X(t,x),t\geq 0,x\in K\cap H_c)$ with $X(0,x) = x$ and a set $K_0$ dense in $K\cap H_c$, such that for all $x \in K_0$ there exists a unique strong solution of equation \eqref{Eq:0.10} given by $\left(\left(X(t,x)\right)_{t\geq0},\eta^x,W\right)$.
\item[ii)] The law of $\left(X(t,x)_{t\geq 0},\eta^x\right)$ is a regular conditional distribution of the law of $\left(\hat{X}_c,\eta\right)$ given $\hat{X}_{c}(0)=x \in K \cap H_c$.
\end{enumerate}
\end{Coro}
\Proof By Theorem \ref{Th:3.1}, we have a stationary strong solution $\hat{X}_{c}$ in $H_c$, such that $W$ and $\hat{X}_{c}(0)$ are independent. Conditioning $\left(\hat{X}_{c},\eta\right)$ on the value of $\hat{X}_{c}(0)=x$, with $c=\overline{x}$, we obtain for $\nu_c$-almost every $x$ a strong solution that we denote $\left(X(t,x),\eta^x\right)$ for all $t\geq0$ and for all $x\in K\cap H_{c}$. This process is the desired process. Indeed, since the support of $\nu_c$ is $K\cap H_c$, we have a strong solution for a dense set $K_0$ in $K \cap H_c$.\\
Notice that all processes $\left(X(t,x)\right)_{t\geq 0}$ with $x \in K_0$ are driven by the same noise $W$ and are continuous with values in $H$. Arguing as in the proof of Lemma \ref{Le:3.1} we see that for all $x,y \in K_0$, for all $t \geq 0$: 
\BD
\|X(t,x)-X(t,y)\|_{-1} \leq \|x-y\|_{-1}.
\ED
Then by density, we obtain a continuous process $\left(X(t,x)\right)_{t\geq 0}$ in $H_c$ for all $x \in K \cap H_c$.
\EProof
We want to prove that for any deterministic initial condition $x \in K\cap H_c$ where $c > 0$, there exists a strong solution of equation \eqref{Eq:0.10}, necessarily unique and that the process $X$ constructed in Corollary \ref{Co:3.1} is a realization of such solution. We have proved this result only for $x$ in a dense set $K_0$, but thanks to the convergence of the transition semigroup $P^{n,c}$, we will be able to conclude. First we prove that the transition semigroup converges on $K \cap H_c$. This result is explained by the following proposition : 
\begin{Prop}\label{Pr:3.3}
Let $c>0$, for all $\phi \in \mathcal{C}_b(H)$ and $ x \in K \cap H_c$:
\BE\label{Eq:3.38}
\lim_{n \rightarrow +\infty} P^{n,c}_t\phi(x) = \E[\phi(X(t,x))]=:P^c_t\phi(x).
\EE
Moreover the Markov process $(X(t,x),t\geq 0, x\in K\cap H_c)$ is strong Feller and its transition semigroup $P^c$ is such that:
\BE\label{Eq:3.39}
|P^c_t\phi(x)-P^c_t\phi(y)|\leq \frac{\|\phi\|_{\infty}}{\sqrt{t}}\|x-y\|_H,\quad \text{ for all }x,y \in K \cap H_c,\text{for all } t>0.
\EE
\end{Prop}
\Proof By proposition \ref{Prop:3.1} $X^n$ is strong Feller on $H_c$ and for all $\phi:H_c\rightarrow\R$ bounded and Borel we have :
\BE\label{Eq:3.40}
|P_t^{n,c}\phi(x) - P_t^{n,c}\phi(y)|\leq \frac{\|\phi\|_{\infty}}{\sqrt{t}}\|x-y\|_{H},\quad \text{ for all } x,y \in K \cap H_c, \text{ for all } t>0. 
\EE
Since $(\nu_{c}^n)_{n\geq 1}$ is tight in $H_{c}$, then there exists an increasive sequence of compact sets $(J^p)_{p\in\N}$ in $H$ such that:
\BE
\lim_{p\rightarrow +\infty} \sup_{n\geq1} \nu_{c}^n(H\setminus J^p) = 0.
\EE
Set $J:= \mathop{\cup}_{p\in\N} J^p\cap K$. Since the support of $\nu_{c}$ is in $K\cap H_{c}$ and $\nu_{c}(J)=1$, then $J$ is dense in $K\cap H_{c}$.
Fix $t>0$, by \eqref{Eq:3.40}, for any $\phi \in \mathcal{C}_b(H)$ :
\BE
\sup_{n \in \N}(\|P^{n,c}_t\phi\|_{\infty} + [P^{n,c}_t\phi]_{Lip(H_c)}) < +\infty.
\EE
Let $(n_j)_{j\in\N}$ be any sequence in $\N$. With a diagonal procedure, by Arzel\`a-Ascoli Theorem, there exists $(n_{j_{l}})_{l\in\N}$ a subsequence and a function $\Theta_{t}:J\rightarrow \R$ such that:
\BE
\lim_{l\rightarrow +\infty} \sup_{x\in J^p}|P^{n_{j_{l}},c}_t\phi(x) - \Theta_{t}(x)|=0, \quad \text{ for all }p\in\N.
\EE
By density, $\Theta_{t}$ can be extended uniquely to a bounded Lipschitz function $\tilde{\Theta}_{t}$ on $K\cap H_{c}$ such that
\BE
\tilde{\Theta}_{t}(x)=\lim_{l\rightarrow + \infty}P^{n_{j_l},c}_t\phi(x),\quad\text{ for all } x \in K\cap H_c.
\EE
Note that the subsequence depends on $t$. Therefore, we have to prove that the limit defines a semigroup and does not depend on the chosen subsequence.\\
By the theorem \ref{Th:3.1}, we have for all $\phi,\psi \in \mathcal{C}_b(H)$ :
\BEAS
\E\left[\psi\left(\hat{X}_{c}(0)\right)\phi\left(\hat{X}_{c}(t)\right)\right] &=&\lim_{l\rightarrow +\infty}\E\left[\psi\left(\hat{X}^{n_{j_{l}}}_{c}(0)\right)\phi\left(\hat{X}^{n_{j_{l}}}_{c}(t)\right)\right]\\
&=& \lim_{l\rightarrow +\infty}\int_H\psi(y)\E\left[\phi\left(\hat{X}^{n_{j_{l}}}_{c}(t)\right)\Big|\hat{X}^{n_{j_{l}}}_{c}(0)=y\right]\nu^{n_{j_{l}}}_c(dy)\\
&=& \lim_{l\rightarrow +\infty}\int_H\psi(y)P^{n_{j_{l}},c}_{t}\phi(y)\nu^{n_{j_{l}}}_c(dy)\\
&=& \int_H\psi(y)\tilde{\Theta}_{t}(y)\nu_c(dy).
\EEAS
Thus, by Corollary \ref{Co:3.1}, we have the following equality:
\BE\label{Eq:3.44}
\E\left[\phi\left(X(t,x)\right)\right]=\tilde{\Theta}_{t}(x), \quad \text{ for }\nu_c \text{-almost every } x.
\EE
Since $\E[\phi(X(t,.))]$ and $\tilde{\Theta}_{t}$ are continuous on $K\cap H_c$, and $\nu_c(K\cap H_c)=1$, the equality \eqref{Eq:3.44} is true for all $x \in K\cap H_c$. Moreover the limit does not depend on the chosen subsequence, and we obtain \eqref{Eq:3.38}. Since the semigroups are equi-Lipschitz, we deduce \eqref{Eq:3.39}.
\EProof

\subsection{Existence of solutions}
We have proved that there exists a continous process $X$ which is a strong solution of equation \eqref{Eq:0.10}  for an $x$ in a dense space. In this section, we prove existence for an initial condition in $K\cap H_c$ with $c>0$. 
\begin{Th}\label{Th:Existence}
Let $\xi$ be a $K$-valued random value with $\overline{\xi}>0$ almost surely and $(\xi,W)$ independent, then there exists a continuous process denoted $(X(t,\xi))_{t\geq 0}$ and a measure $\eta^\xi$ such that: \begin{enumerate}
\item[(a)] $\left(\left(X(t,\xi)\right)_{t\geq0},\eta^\xi,W\right)$is the unique strong solution of \eqref{Eq:0.10} with $X(0,\xi) = \xi$ almost surely.
\item[(b)] The Markov process $(X(t,x),t\geq 0, x\in K\cap H_c)$ is continous and has $P^c$ for transition semigroup which is strong Feller on $H_c$.
\item[(c)] For all $c >0$, $x \in K \cap H_c$ and $0=t_0<t_1<\cdots <t_m$, $(X(t_i,x) , i =1, \dots, n)$ is the limit in distribution of $(X^n(t_i,x))_{i =1, \dots, m}$.
\item[(d)] If $\xi$ has distribution $\nu_c$ with $c>0$, then $\left(X(t,\xi)\right)_{t\geq0}$ is equal in distribution to $\left(\hat{X}_{c}(t)\right)_{t\geq 0}$.
\end{enumerate}
\end{Th}
\Proof By Corollary \ref{Co:3.1} we have a process $(X(t,x),t\geq0,x\in K \cap H_c)$, such that for all $x$ in a set $K_0$ dense in $K \cap H_c$ we have a strong solution $\left(\left(X(t,x)\right)_{t\geq0},\eta^x,W\right)$ of \eqref{Eq:0.10} with initial condition $x$. By proposition \ref{Pr:3.3}, we have that the Markov process $X$ has transition semigroup $P_c$ on $H_c$.\\
The strong Feller property of $P^c$ implies that for all $x\in K \cap H_c$ and $s>0$ the law of $X(s,x)$ is absolutely continous with respect to the invariant measure $\nu_c$. Indeed, if $\nu_{c}(\Gamma)=0$, then $\nu_c(P_s^c(\1_{\Gamma})) = \nu_c(\Gamma)=0$. So $P_s^c(\1_{\Gamma})(x)=0$ for $\nu_c$-almost every $x$ and by continuity for all $x \in K \cap H_c$.\\
Therefore almost surely $X(s,x) \in K_0$ for all $s>0$ and $x \in K \cap H_c$. Fix $s>0$, denote for all $\theta \in [0,1]$:
\BD
\tilde{X} := t \mapsto X(t+s,x),
\tilde{W}(\cdot,\theta) := t \mapsto W(t+s,\theta)-W(s,\theta)),
\ED
and the measure $\tilde{\eta}^x$ such that for all $T>0$, and for all $h \in \mathcal{C}(O_{T})$:
\BD
\big\langle h, \tilde{\eta}^x\big\rangle_{O_{T}} := \int_{O_{s}^{T+s}} h(t-s,\theta)\eta^x(dt,d\theta)
\ED
So we have a process $\tilde{X} \in \mathcal{C}([0,T];H) \cap \mathcal{C}(O_{T})$ and a mesure $\tilde{\eta}^x$ on $O_{T}$ which is finite on $[\delta,T]\times[0,1]$ for all $\delta \geq0$, such that $\left((\tilde{X}(t,x))_{t\geq0},\tilde{\eta}^x,\tilde{W}\right)$ is a strong solution of \eqref{Eq:0.10} with initial condition $X(s,x)$. By continuity $X(s,x)\rightarrow x$ in $H$ as $s\rightarrow0$, so $\left((X(t,x))_{t\geq0},\eta^x,W\right)$ is a strong solution of \eqref{Eq:0.10} with initial condition $x$ in the sense of the definition \ref{De:2}.\\
Thanks to the previous results, (b), (c) and (d) are obvious.
\EProof

\newpage

\section{Reflection and Revuz measures}\label{S:5}
We have proved the existence of solution to \eqref{Eq:0.10} with a reflection measure. In \cite{MR1959795}, L. Zambotti uses an integration by parts formula to prove that, in some cases, the reflection measure vanishes. Moreover, L. Zambotti proves that, in some other cases, the reflection measure does not vanish. He uses the theory of the Continuous Additive Functionnals described in \cite{MR1303354}. We adapt his arguments and prove similar results for our case.
\subsection{Integration by parts formula}\label{S:5.1}
For all $\phi \in \mathcal{C}_{b}^1(H_{c})$ we denote by $\partial_{h}\phi$ the directional derivative of $\phi$ along $h\in H$:
\BD
\partial_{h}\phi :x \mapsto \lim_{t\rightarrow0}\frac{1}{t}(\phi(x+th)-\phi(x)), \quad x \in H.
\ED
For all $\phi \in \mathcal{C}_{b}^1(H)$, we have:
\BD
\langle\nabla \phi (x),h\rangle = \partial_{h}\phi(x).
\ED
 We have the following classical result (see \cite{MR2349572} for details):
We denote by $(M,\hat{M})$ two indepedant copies of the standard Brownian meander (see  \cite{MR1083357} and \cite{MR726906}), and we set for all $r \in (0,1)$:
\BE
\mathcal{U}_{r}(\theta):=\left\{
\begin{array}{lr}
\sqrt{r}M\left(\frac{r-\theta}{r}\right),& \theta \in [0,r],\\
\sqrt{1-r}\hat{M}\left(\frac{\theta-r}{1-r}\right),& \theta \in ]r,1].
\end{array}\right.
\EE
The starting point is the Theorem 7.1 in \cite{MR2349572} where the following formula has been proved for a process $Y$ whose the law is $\mu$.
\begin{Th}
For all $\Phi$ in $\mathcal{C}^1_{b}(H,\R)$ and $h\in D(A)$:
\BEA
\E\left[\partial_{h}\Phi(Y)\1_{Y\in K}\right]&=&-\E\left[\left(\langle Y,Ah\rangle - \overline{Y}\cdot\overline{h}\right)\Phi(Y)\1_{Y\in K}\right]\\
&&-\int_{0}^{1}h(r)\frac{1}{\sqrt{2\pi^3r(1-r)}}\E\left[\Phi(\mathcal{U}_{r})e^{-(1/2)(\overline{\mathcal{U}}_{r})^2}\right]dr.\nonumber
\EEA
\end{Th}
We denote by $p_{\overline{\mathcal{U}}_{r}}:\R^+\rightarrow[0,1]$ the continuous version of the density of $\overline{\mathcal{U}}_{r}$. By conditioning on $\overline{Y}=c$, we obtain:
\BEA\label{Eq:4.14}
\E\left[\partial_{\Pi h}\Phi(Y_{c})\1_{Y_{c}\in K}\right]&=&-\E\left[\langle Y_{c},Ah\rangle\Phi(Y_{c})\1_{Y_{c}\in K}\right]\\
&&-\int_{0}^{1}\Pi h(r)\frac{p_{\overline{\mathcal{U}}_{r}}(c)}{\pi\sqrt{r(1-r)}}\E\left[\Phi(\mathcal{U}_{r})|\overline{\mathcal{U}}_{r}=c\right]dr,\nonumber
\EEA
where $Y_{c}$ has been defined in the section \ref{S:1}. Moreover, notice that we have the following classical and easy to prove integration by parts formula for the measures $(\nu_{c}^n)_{n\in\N}$. For all $\Phi$ in $\mathcal{C}^1_{b}(H)$ and $h\in D(A)$:
\BE\label{Eq:4.15}
\int_{H}\partial_{\Pi h}\Phi\ d\nu_{c}^n = -\int_{H}\langle x, Ah \rangle \Phi(x) \nu_{c}^n(dx) - \int_{0}^1 \Pi h(r) \int_{H}\Phi(x) f^n(x(r))\nu_{c}^n(dx) dr.
\EE
We define $\gamma^n : x \mapsto \frac{1}{Z^n_c}\exp(-U^n(x))$ for all $x\in H$, where $Z^n_{c}$ is the constant of normalization defined in \eqref{Eq:2.10}. Then $\gamma^n \in \mathcal{C}^1_{b}(H)$ and for all $x,h \in K$:
\BE
\langle \nabla \gamma^n(x),h\rangle =  \gamma^n(x)\langle \nabla \log \gamma^n(x),h\rangle= \gamma^n(x)\int_{0}^1 h(\theta) f^n(x(\theta))d\theta.
\EE
Let $\phi$ be  in $\mathcal{C}^1_{b}(H)$. We use \eqref{Eq:4.14}, with $\Phi = \phi\cdot \gamma^n$. So we obtain:
\BEA\label{Eq:4.17}
\int_{H}\partial_{\Pi h}(\phi\cdot \gamma^n)\ d\mu_{c} &=& -\int_{H}\langle x, Ah \rangle \phi(x)\gamma^n(x) \1_{x\in K}\mu_{c}(dx)\\
&& - \int_{0}^{1}\Pi h(r)\frac{p_{\overline{\mathcal{U}}_{r}}(c)}{\pi\sqrt{r(1-r)}}\E\left[\phi(\mathcal{U}_{r})\gamma^n(\mathcal{U}_{r})|\overline{\mathcal{U}}_{r}=c\right]dr.\nonumber
\EEA
We compute the derivative of the product, and obtain:
\BEA\label{Eq:4.18}
\int_{H}(\partial_{\Pi h}\phi)\cdot \gamma^n\ d\mu_{c} &=& -\int_{H}\Big(\langle x, Ah \rangle +\langle \nabla \log \gamma^n(x),\Pi h\rangle\Big)\phi(x)\gamma^n(x) \1_{x\in K}\mu_{c}(dx)\\
&& - \int_{0}^{1}\Pi h(r)\frac{p_{\overline{\mathcal{U}}_{r}}(c)}{\pi\sqrt{r(1-r)}}\E\left[\phi(\mathcal{U}_{r})\gamma^n(\mathcal{U}_{r})|\overline{\mathcal{U}}_{r}=c\right]dr.\nonumber
\EEA
We want to let $n$ go to infinity. We have to study the convergence of all the terms. \\
By section \ref{S:3.2}, the left-hand side  converges to:
\BD
\int_{H}(\partial_{\Pi h}\phi )d\nu_{c}.
\ED
Denote now by $I_{r}^n$ the following term:
\BD
I_{r}^n:=\frac{p_{\overline{\mathcal{U}}_{r}}(c)}{\pi\sqrt{r(1-r)}}\E\left[\phi(\mathcal{U}_{r})\gamma^n(\mathcal{U}_{r})\big|\overline{\mathcal{U}}_{r}=c\right].
\ED
Since $Z^n_{c}$ converges, there exists $C$ such that for all $r \in (0,1)$:
\BD
|I_{r}^n| \leq C\frac{p_{\overline{\mathcal{U}}_{r}(c)}}{\sqrt{r(1-r)}}\|\phi\|_{\infty}J_{r}^n
\ED
where $J_{r}^n$ is defined by:
\BD
J_{r}^n:=\E\left[\exp\left(-\int_{0}^1F^n(\mathcal{U}_{r}(\theta)) d\theta\right)\right].
\ED
In the logarithmic case and in the negative $\alpha$-power case, as in section \ref{S:3.2} and by dominated convergence, we have for all $r\in (0,1)$:
\BE
\lim_{n\rightarrow +\infty}J_{r}^n = \E\left[\exp\left(-\int_{0}^1F(\mathcal{U}_{r}(\theta)) d\theta\right)\right].
\EE
Therefore, in the logarithmic case and in the negative $\alpha$-power case for $\alpha>1$, since $|J^{n}_{r}|<1$, by dominated convergence, the last term in \eqref{Eq:4.18} has a limit when $n$ grows to infinity.\\
In the negative $\alpha$-power case for $\alpha\leq1$, since 
\BD
|J^{n}_{r}|\leq\E\left[\exp\left(-\int_{0}^1F^1_{\alpha}(\mathcal{U}_{r}(\theta)) d\theta\right)\right],
\ED
by dominated convergence, the last term in \eqref{Eq:4.18} has a limit when $n$ grows to infinity.\\
Moreover, if $\alpha\geq 3$, by the law of the iterated logarithm, almost surely and for all $r\in (0,1)$:
\BD
\int_{0}^1\frac{d\theta}{(\mathcal{U}_{r}(\theta))^{\alpha-1}} =+\infty.
\ED
Thus, in this case
\BE
\lim_{n\rightarrow +\infty}J_{r}^n = 0,
\EE
and, by dominated convergence, the last term in \eqref{Eq:4.18} converges to $0$.\\
Now we use the representation described in \cite{MR2349572} in order to prove the convergence of the first term in the right-hand side of \eqref{Eq:4.18}. Denote by $S^n$ the following
\BE
\begin{array}{rcl}
S^n &:=&-\int_{H}\left(\langle x, Ah \rangle +\langle \nabla \log \gamma^n(x),\Pi h\rangle\right)\phi(x)\gamma^n(x) \1_{x\in K}\mu_{c}(dx)\\
\\
&=& -\E \left[ \left(\langle Y_{c}, Ah \rangle +\langle \nabla \log \gamma^n(Y_{c}),\Pi h\rangle\right)\phi(Y_{c})\gamma^n(Y_{c})\1_{Y_{c}\in K} \right]
\end{array}
\EE
We use the following Theorem whose proof is in Appendix A in \cite{MR2349572}.
\begin{Th}
For all $\Psi : \mathcal{C}([0,1]) \rightarrow \R$ bounded and Borel
\BE
\begin{array}{rcl}
\E[\Psi(Y)] &=& \int_{\R} \E\left[\Psi(y+\B)\sqrt{\frac{4}{3}}\exp\left(-\frac{1}{2}(y+\overline{\B})^2 + \frac{3}{8}y^2\right)\right]\frac{\exp\left(-\frac{3}{8}y^2\right)}{\sqrt{2\pi 4/3}}dy\\
&=& \int_{\R}\E\left[\Psi(y+\B)\rho(y+\B)\right]dy,
\end{array}
\EE
where we set $\rho : \mathcal{C}([0,1]) \rightarrow \R$,
\BD
\rho(u) :=\frac{1}{\sqrt{2\pi}}\exp\left(-\frac{1}{2}\overline{u}^2\right).
\ED
\end{Th}
Thanks to this Theorem, we can write:
\BD
\begin{array}{rcl}
S^n &=&- \int_{\R}\E \Big[ \left(\left\langle y+\B, Ah \right\rangle +\left\langle \nabla \log \gamma^n(y+\B),\Pi h\right\rangle\right)\\
\\
&&\times\phi(y+\B)\gamma^n(y+\B) \rho(y+\B)\1_{y+\B\in K}\Big| \overline{\B}=c-y\Big]dy
\end{array}
\ED
We set $\mathcal{V}_{r} = -\sqrt{r}M(1) + \mathcal{U}_{r}$. Notice that $\mathcal{V}_{r}$ is $0$ at time $0$, then run backwards the path of $M$ on $[0,r]$ and then runs the path of $\hat{M}$ on $]r,1]$. Almost surely since $M>0$ on $]0,1]$, then $\mathcal{V}_{r}$ attains the minimum $-\sqrt{r}M(1)$ only at time $r$. Let $(\tau, M, \hat{M})$ be an independent triple, such that $\tau$ has the arcsine law, then $\mathcal{V}_{\tau}$ has the same law as $\B$ (see \cite{MR726906}). We can write:
\BD
\begin{array}{rcl}
S^n &=&-\int_{0}^1\frac{1}{\pi\sqrt{r(1-r)}} \int_{\R}\E \Big[ \left(\langle y+\mathcal{V}_{r}, Ah \rangle +\left\langle f\left(\frac{1}{n}+y+\mathcal{V}_{r}\right),\Pi h\right\rangle\right)\\
 \\
&&\times\phi(y+\mathcal{V}_{r})\gamma\left(\frac{1}{n}+y+\mathcal{V}_{r}\right) \rho(y+\mathcal{V}_{r})\1_{y+\mathcal{V}_{r}\in K}\Big|  \overline{\mathcal{V}}_{r}=c-y\Big]dy\ dr\\
\\
&=&-\int_{0}^1\frac{1}{\pi\sqrt{r(1-r)}} \int_{\R}\E \Big[ \left(\langle z-\frac{1}{n}+\mathcal{V}_{r}, Ah \rangle +\left\langle f\left(z+\mathcal{V}_{r}\right),\Pi h\right\rangle\right)\\
 \\
&&\times\phi(z-\frac{1}{n}+\mathcal{V}_{r})\gamma\left(z+\mathcal{V}_{r}\right) \rho(z-\frac{1}{n}+\mathcal{V}_{r})\1_{z-\frac{1}{n}+\mathcal{V}_{r}\in K}\Big|  \overline{\mathcal{V}}_{r}=c-z+\frac{1}{n}\Big]dz\ dr.
\end{array}
\ED
Now we use the proposition \ref{Prop:5.4} which is stated in the next section \ref{S:5.3}. Thus, we can used Fatou Lemma to prove that for all $h \in D(A)$:
\BD
\frac{1}{\pi\sqrt{r(1-r)}}\left[ \left\langle f(z+\mathcal{V}_{r}),\Pi h\right\rangle\|\phi\|_{\infty}\gamma\left(z+\mathcal{V}_{r}\right)\1_{z+\mathcal{V}_{r}\in K}\right]
\ED
is integrable on $\Omega\times\R \times [0,1]$.
Thus, we can used the dominated convergence Theorem to see:
\BE
\begin{array}{rcl}
\lim_{n\rightarrow+\infty}S^n  &=&-\int_{0}^1\frac{1}{\pi\sqrt{r(1-r)}} \int_{\R}\E \Big[ \Big(\langle z+\mathcal{V}_{r}, Ah \rangle +\left\langle f\left(z+\mathcal{V}_{r}\Big),\Pi h\right\rangle\right)\\
 \\
&&\times\phi(z+\mathcal{V}_{r})\gamma\left(z+\mathcal{V}_{r}\right) \rho(z+\mathcal{V}_{r})\1_{z+\mathcal{V}_{r}\in K}\Big|  \overline{\mathcal{V}}_{r}=c-z\Big]dz\ dr\\
\\
&=& - \int_{\R}\E \Big[ \Big(\langle z+\mathcal{V}_{\tau}, Ah \rangle +\left\langle f\left(z+\mathcal{V}_{\tau}\right),\Pi h\right\rangle\Big)\\
 \\
&&\times\phi(z+\mathcal{V}_{\tau})\gamma\left(z+\mathcal{V}_{\tau}\right) \rho(z+\mathcal{V}_{\tau})\1_{z+\mathcal{V}_{\tau}\in K}\Big|  \overline{\mathcal{V}}_{\tau}=c-z\Big]dz\\
\\
&=& - \E \left[ \Big(\langle Y, Ah \rangle +\left\langle f(Y),\Pi h\right\rangle\Big)\phi(Y)\gamma\left(Y\right) \1_{Y\in K}\Big|  \overline{Y}=c\right]\\
\\
&=& - \E \left[ \Big(\langle Y_{c}, Ah \rangle +\left\langle f(Y),\Pi h\right\rangle\Big)\phi(Y_{c})\gamma\left(Y_{c}\right) \1_{Y_{c}\in K}\right]\\
\\
&=& - \int_{H} \Big(\langle x, Ah \rangle +\left\langle f(x),\Pi h\right\rangle\Big)\phi(x)\nu_{c}(dx)
\end{array}
\EE
For all $r\in (0,1)$, denote $\Sigma_{r}^c$ such that:
\BE
\Sigma_{r}^c(d\omega):=\frac{1}{\mu_{c}(K)}\frac{p_{\overline{\mathcal{U}}_{r}}(c)}{\pi\sqrt{r(1-r)}}\PP\left(\mathcal{U}_{r}\in d\omega|\overline{\mathcal{U}}_{r}=c\right),
\EE
thus we have the following Theorem:
\begin{Th}
For all $\phi$ in $\mathcal{C}^1_{b}(H)$ and $h\in D(A)$:
\BEA\label{Eq:5.14}
\int_{H}\partial_{\Pi h}\phi(x) \1_{x\in K}\nu_{c}(dx)&=&- \int_{H} \Big(\langle x, Ah \rangle +\left\langle f(x),\Pi h\right\rangle\Big)\phi(x)\nu_{c}(dx)\\
&&-\int_{0}^{1}\Pi h(r)\int\phi\gamma d\Sigma_{r}^c\ dr.\nonumber
\EEA
Moreover, for $\alpha\geq3$, the last term vanishes.
\end{Th}

\subsection{Dirichlet forms}
We now describe the Dirichlet Forms and the resolvent associated to $\hat{X}^n_{c}$, in order to obtain the Dirichlet Forms and the resolvent associated to $\hat{X}_{c}$. The first result is the following description of the generator of $Z$. Let $\psi_{h}: x \mapsto \exp(i(x,h)_{-1})$ for $x\in H_{c}$ and $h \in D(A^2)$, then the generator of $Z$ is such that
\BD
L\psi_{h}(x):= \frac{d}{dt}\E[\psi_{h}(Z(t,x))]\Big|_{t=0} = -\frac{1}{2}\psi_{h}(x)\Big(i(A^2h,x)_{-1}+\|\Pi h\|_{-1}^2\Big)
\ED
We define for all $\phi \in \mathcal{C}_{b}(H_{c})$ the resolvent of $\hat{X}^n_{c}$ on $H_{c}$:
\BE
R_{\lambda}^{n,c}\phi(x) := \int_{0}^{\infty} e^{-\lambda t}\E\left[\phi(\hat{X}^n_{c}(t,x))\right]dt, \quad x \in H_{c}, \lambda >0.
\EE
We define $Exp_{A}(H_{c})\subset \mathcal{C}_{b}(H_{c})$ as the linear span of $\{cos((h,\cdot));sin((h,\cdot)), h \in D(A^2)\}$. Then we define the symmetric bilinear form:
\BE
\mathcal{E}^{n,c}(\phi,\psi) := \frac{1}{2}\int_{H} \langle -A \nabla\phi, \nabla\psi\rangle d\nu_{c}^n, \quad\text{ for all } \phi, \psi \in Exp_{A}(H).
\EE
The following result is standard.
\begin{Prop}
$(\mathcal{E}^{n,c}, Exp_{A}(H_{c}))$ is closable in $L^2(\nu_{c}^n)$: we denote by $(\mathcal{E}^{n,c},D(\mathcal{E}^{n,c}))$ the closure.\\
$(R_{\lambda}^{n,c})_{\lambda>0}$ is the resolvent associated with $\mathcal{E}^{n,c}$, that is, for all $\lambda >0$ and $\psi \in L^2(\nu_{c}^n)$, $R_{\lambda}^{n,c}\psi \in D(\mathcal{E}^{n,c})$ and:
\BE
\lambda \int_{H} R_{\lambda}^{n,c}\psi \phi\ d\nu_{c}^n + \mathcal{E}^{n,c}(R_{\lambda}^{n,c}\psi,\phi)=\int_{H}\psi\phi\ d\nu_{c}^n, \quad\text{ for all } \phi \in D(\mathcal{E}^{n,c}).
\EE
\end{Prop}
Let $\psi_{h}: x \mapsto \exp(i(x,h))$ for $x\in H_{c}$ and $h \in D(A^2)$. By It\^o formula
\BE
L^n\psi_{h}(x):= \frac{d}{dt}\E[\psi_{h}(\hat{X}^n_{c}(t,x))]\Big|_{t=0} = L\psi_{h}(x) + \frac{i}{2}\langle f^n(x), \Pi h \rangle \psi_{h}(x).
\EE
After an easy computation, we have $(L^n, Exp_{A}(H_{c}))$ is symmetric in $L^2(\nu_{c}^n)$ and:
\BE
\int_{H} L^n\phi\psi\ d\nu_{c}^n = -\frac{1}{2}\int_{H}\langle-A\nabla\phi, \nabla\psi\rangle d\nu_{c}^n, \quad \text{ for all } \phi, \psi \in Exp_{A}(H_{c}).
\EE
Moreover we define for all $\phi \in \mathcal{C}_{b}(H_{c})$ the resolvent of $\hat{X}_{c}$ on $K \cap H_{c}$:
\BE
R_{\lambda}^{c}\phi(x) := \int_{0}^{\infty} e^{-\lambda t}\E\left[\phi\left(\hat{X}_{c}(t,x)\right)\right]dt, \quad x \in K \cap H_{c}, \lambda >0.
\EE
We also define the symmetric bilinear form:
\BE
\mathcal{E}^{c}(\phi, \psi) := \frac{1}{2}\int_{H} \langle-A \nabla\phi, \nabla\psi\rangle d\nu_{c}, \quad\text{ for all } \phi, \psi \in\mathcal{C}^1_{b}(H).
\EE
Proceeding as in the proposition 8.1 in \cite{MR2349572}, we can prove that for all $\phi,\psi \in \mathcal{C}^1_{b}(H_{c})$, $\mathcal{E}^{n,c}(\phi, \psi) \rightarrow \mathcal{E}^{c}(\phi, \psi)$ and $R_{\lambda}^{n,c}\phi \rightarrow R_{\lambda}^{c}\phi$ uniformly as $n$ grows to infinity. 
Let $\psi\in\mathcal{C}_{b}(H_{c})$, we can write for all $h\in D(A^2)$:
\BE
\begin{array}{rcl}
\int_{H} \psi \phi _{h}d\nu_{c}&=& \lim_{n\rightarrow +\infty} \int_{H} \psi \phi_{h} d\nu_{c}^n = \lim_{n\rightarrow +\infty} \int_{H} R_{\lambda}^{n,c}\psi(\lambda \phi_{h} - L^n \phi_{h}) d\nu_{c}^n\\
\\
&=&\int_{H} R_{\lambda}^{c}\psi(\lambda \phi_{h} - L \phi_{h}) d\nu_{c}-\frac{i}{2}\lim_{n\rightarrow +\infty} \int_{H} R_{\lambda}^{n,c}\psi(x)\phi_{h}(x)\langle f^n(x),\Pi h\rangle d\nu_{c}^n.
\end{array}
\EE
Then, with the proposition \ref{Prop:5.4} below:
\BE
\begin{array}{rcl}
\int_{H} \psi \phi _{h}d\nu_{c}&=&\int_{H} R_{\lambda}^{c}\psi(\lambda \phi_{h} - L \phi_{h}) d\nu_{c}-\frac{i}{2} \int_{H} R_{\lambda}^{c}\psi(x)\phi_{h}(x)\langle f(x),\Pi h\rangle d\nu_{c}\\
\\
&&-\frac{i}{2}\int_{0}^1\Pi h(r) \int_{H} R_{\lambda}^{c}\psi \phi_{h} \gamma\ d\Sigma_{r}^c\ dr\\
\\
&=&\lambda\int_{H} R_{\lambda}^{c}\psi\phi_{h} d\nu_{c} - \frac{i}{2} \int_{H} R_{\lambda}^{c}\psi \phi_{h} \langle A h, x\rangle d\nu_{c} + \frac{1}{2} \int_{H} R_{\lambda}^{c}\psi \phi_{h} \|\Pi h\|_{-1}^2 d\nu_{c}\\
\\
&&-\frac{i}{2} \int_{H} R_{\lambda}^{c}\psi(x)\phi_{h}(x)\langle f(x),\Pi h\rangle d\nu_{c}-\frac{i}{2}\int_{0}^1\Pi h(r) \int_{H} R_{\lambda}^{c}\psi \phi_{h} \gamma\  d\Sigma_{r}^c\ dr.
\end{array}
\EE
Thanks to the integration by parts formula applied to $R_{\lambda}^{c}\psi\phi_{h}$ , we have:
\BE
\int_{H} \psi \phi _{h}d\nu_{c}=\lambda\int_{H} R_{\lambda}^{c}\psi\phi_{h} d\nu_{c} + \mathcal{E}^{c}(R_{\lambda}^{c}\psi, \phi_{h}).
\EE
By linearity and by density, we obtain for all $\lambda >0$ and $\psi \in \mathcal{C}_{b}(H)$:
\BE
\lambda \int_{H} R_{\lambda}^{c}\psi \phi\ d\nu_{c} + \mathcal{E}^{c}(R_{\lambda}^{c}\psi,\phi)=\int_{H}\psi\phi\ d\nu_{c}, \quad \text{ for all } \phi \in \mathcal{D},
\EE
where we denote $\mathcal{D}:=\{R_{\lambda}^{c}\phi, \phi \in \mathcal{C}_{b}(H_{c}),\lambda >0\}$. We use classical results from \cite{MR1214375}, and obtain the following proposition:
\begin{Prop} Let $c>0$.\\
i) $(\mathcal{E}^{c}, Exp_{A}(H_{c}))$ is closable in $L^2(\nu_{c})$: we denote by $(\mathcal{E}^{c},D(\mathcal{E}^{c}))$ the closure.\\
ii) $(\mathcal{E}^{c},D(\mathcal{E}^{c}))$ is a symmetric Dirichlet form such that $Lip(H_{c}) \subset D(\mathcal{E}^c)$ and $\mathcal{E}^c(\phi,\phi) \leq |\phi|_{Lip(H_{c})}^2$.
iii) $(R_{\lambda}^{c})_{\lambda>0}$ is the resolvent associated with $\mathcal{E}^{c}$, that is, for all $\lambda >0$ and $\psi \in L^2(\nu_{c})$, $R_{\lambda}^{c}\psi \in D(\mathcal{E}^{c})$ and:
\BE
\lambda \int_{H} R_{\lambda}^{c}\psi \phi\ d\nu_{c} + \mathcal{E}^{c}(R_{\lambda}^{c}\psi,\phi)=\int_{H}\psi\phi\ d\nu_{c}, \quad\text{ for all } \phi \in D(\mathcal{E}^{c}).
\EE
iv) $(P_{t}^c)_{t\geq0}$ is the semigroup associated with $(\mathcal{E}^{c},D(\mathcal{E}^{c}))$.
\end{Prop}

\subsection{Total mass of the reflection measure}\label{S:5.3}
We now state and prove the proposition \ref{Prop:5.4} used above.
\begin{Prop}\label{Prop:5.4} For all $\phi \in \mathcal{C}_{b}(H_{c})$, for all $h\in D(A)$:
\BE\label{Eq:5.27}
\int_{H}\langle f^n(x),h\rangle \phi(x)\gamma^n(x)\mu_{c}(dx) 
\EE
has a limit when $n$ grows to infinity.\\
Moreover for all $0<\delta\leq s\leq t \leq T$, $\E\left[\eta(O_{s,t})\right]<+\infty$.
\end{Prop}
\Proof 
Denote $\sigma^n_{r,c}$ the measure such that for all $r \in [0,1]$, for all $c>0$:
\BD
\sigma_{r,c}^n(dx):= f^n\left(x(r)\right) \gamma^n(x)\mu_{c}(dx) 
\ED
It suffices to prove that:
\BE
\limsup_{n\rightarrow+\infty}\left|\int_{0}^1\int_{H} d\sigma^n_{r,c}\ dr\right|<+\infty.
\EE
By symetry, it suffices to prove convergence of 
\BE\label{Eq:5.28}
\left|\int_{0}^{1/2}\int_{H} d\sigma^n_{r,c}\ dr\right|<+\infty.
\EE
The idea is to study an integration by parts formula for the law of $Y_{c}$ on the path space
\BD
\tilde{K}:=\{h\in \mathcal{C}([0,1]),h(\theta)\geq0 \text{ for all } \theta \in [0,1/2]\}.
\ED
The crucial tool is that, on this space, the processes that we consider have no more fixed mean, and we can have an integration by parts formula without the constraint of zero mean.
We set for all $r \in (0,1/2)$:
\BD
\mathcal{T}_{r}(\theta):=\left\{
\begin{array}{lr}
\sqrt{r}M\left(\frac{r-\theta}{r}\right),& \theta \in [0,r],\\
\sqrt{\left(\frac{1}{2}-r\right)}\hat{M}\left(\frac{\theta-r}{\frac{1}{2}-r}\right),& \theta \in ]r,1/2].
\end{array}\right.
\ED
Moreover we set 
\BD
\chi : \theta\mapsto \1_{[0,1/2]}(\theta),
\ED
and for $u\in\mathcal{C}([0,1/2])$
\BD
m(u) :=\int_{0}^{1/2}(u(\theta)+u(1/2))d\theta.
\ED
The starting point is the Lemma B.1 in \cite{MR2349572} where the following formulae have been proved.
\begin{Le}
For all $\Psi:\mathcal{C}([0,1/2])\rightarrow \R$ bounded and Borel:
\BEA
\E\left[\Psi(Y_{c})\right]&=&\sqrt{32}\ \E\left[\Psi(b+\B)\exp\left(-12(m(b+\B)-c)^2+\frac{3}{8}b^2\right)\right]\\
&=& \int_{\R}\E\left[\Psi(y+\B)\tilde{\rho}(y+\B)\right]dy,
\EEA
where we set $\tilde{\rho} : \mathcal{C}([0,1/2]) \rightarrow \R$,
\BD
\tilde{\rho}(u) :=\sqrt{\frac{12}{\pi}}\exp\left(-12(m(u)-c)^2\right).
\ED
Moreover, for all $c>0$ and $\Phi\in\mathcal{C}^1_{b}(L^2(0,1/2))$:
\BEA\label{Eq:5.29}
\E\left[\partial_{\chi}\Phi(Y_{c})\1_{Y_{c}\in \tilde{K}}\right]&=&\E\left[24\left(m(Y_{c})-c\right)\Phi(Y_{c})\1_{Y_{c}\in \tilde{K}}\right]\\
&&-\int_{0}^{1/2}\sqrt{\frac{12}{\pi^3\sqrt{r(1/2-r)}}}\E\left[\Phi(\mathcal{T}_{r})e^{-12(m(\mathcal{T}_{r})-c)^2}\right]dr.\nonumber
\EEA
\end{Le}
We have writen $\Phi(Y_{c})$ for $\Phi(Y_{c}|_{[0,1/2]})$ with a slight abuse of notation.
We set now for $n\geq 1$, $r\in (0,1/2)$:
\BD
\tilde{U}^n(x):=\int_{0}^{1/2}F^n(x(\theta))d\theta,\quad x \in L^2(0,1).
\ED
We define $\tilde{\gamma}^n : x \mapsto \exp(-\tilde{U}^n(x))$ for all $x\in H$.
Then $\tilde{\gamma}^n \in \mathcal{C}^1_{b}(L^2(0,1/2))$ and for all $x,h \in \tilde{K}$:
\BE
\langle \nabla \tilde{\gamma}^n(x),h\rangle =  \tilde{\gamma}^n(x)\langle \nabla \log \tilde{\gamma}^n(x),h\rangle= \tilde{\gamma}^n(x)\int_{0}^{1/2} h(\theta) f^n(x(\theta))d\theta.
\EE
Moreover we define for $n\geq 1$, $r\in (0,1/2)$ and $\Psi \in\mathcal{C}^1_{b}(L^2(0,1/2))$:
\BD
\tilde{\Sigma}^{n,c}_{r}(\Psi):=\sqrt{\frac{12}{\pi^3\sqrt{r(1/2-r)}}}\E\left[\Psi(\mathcal{T}_{r})\tilde{\gamma}^n(\mathcal{T}_{r})e^{-12(m(\mathcal{T}_{r})-c)^2}\right].
\ED
Let $\phi$ be  in $\mathcal{C}^1_{b}(L^2(0,1/2))$. We use \eqref{Eq:5.29}, with $\Phi = \phi\cdot \gamma^n$. So we obtain:
\BE\label{Eq:5.31}
\int_{H}\partial_{\chi}(\phi\cdot \tilde{\gamma}^n)\1_{\tilde{K}} d\mu_{c} = \E\left[24\left(m(Y_{c})-c\right)\phi(Y_{c})\tilde{\gamma}^n(Y_{c})\1_{Y_{c}\in \tilde{K}}\right] -\int_{0}^{1/2}\tilde{\Sigma}^{n,c}_{r}(\phi)dr.
\EE
We compute the derivative of the product, and take $\phi \equiv 1$, then we obtain:
\BE\label{Eq:5.33}
\E\left[\langle \nabla \log \tilde{\gamma}^n(x),\chi\rangle\tilde{\gamma}^n(Y_{c})\1_{Y_{c}\in \tilde{K}}\right] = \E\left[24\left(m(Y_{c})-c\right)\tilde{\gamma}^n(Y_{c})\1_{Y_{c}\in \tilde{K}}\right]-\int_{0}^{1/2}\tilde{\Sigma}^{n,c}_{r}(1)dr.
\EE
Define now for $n\geq 1$, $r\in (0,1/2)$:
\BD
\tilde{U}(x):=\int_{0}^{1/2}F(x(\theta))d\theta,\quad x \in L^2(0,1).
\ED
We also define $\tilde{\gamma} : x \mapsto \exp(-\tilde{U}(x))$ for all $x\in H$. Moreover we define for $n\geq 1$, $r\in (0,1/2)$ and $\Psi \in\mathcal{C}^1_{b}(L^2(0,1/2))$:
\BD
\tilde{\Sigma}^{c}_{r}(\Psi):=\sqrt{\frac{12}{\pi^3\sqrt{r(1/2-r)}}}\E\left[\Psi(\mathcal{T}_{r})\tilde{\gamma}(\mathcal{T}_{r})e^{-12(m(\mathcal{T}_{r})-c)^2}\right].
\ED
Finally, we denote $\tilde{\sigma}^n_{r,c}$ the measure such that for all $r \in [0,1]$, for all $c>0$:
\BD
\tilde{\sigma}_{r,c}^n(dx):= f^n\left(x(r)\right) \tilde{\gamma}^n(x)\mu_{c}(dx) 
\ED
We easily prove the following result:
\begin{Le} For all $c>0$:
\BD
\lim_{n\rightarrow+\infty}\int_{0}^{1/2}\int_{H} d\tilde{\sigma}^n_{r,c}=\E\left[24\left(m(Y_{c})-c\right)\tilde{\gamma}(Y_{c})\1_{Y_{c}\in \tilde{K}}\right]-\int_{0}^{1/2}\tilde{\Sigma}^{c}_{r}(1)dr.
\ED
Moreover, for $\alpha \geq 3$ the last term vanishes.
\end{Le}
We set now for $n\geq 1$:
\BD
\tilde{U'}^n(x):=\int_{1/2}^{1}F^n(x(\theta))d\theta=U^n(x) - \tilde{U}^n(x),\quad x \in L^2(0,1).
\ED
We also define $\tilde{\gamma}^{,n} : x \mapsto \exp(-\tilde{U'}^n(x))$ for all $x\in H$.\\
We notice now that we can compute explicitly the conditional distribution of $Y_{c}$ given $(Y_{c}(\theta),\theta \in [0,1/2])$. Indeed, we have for all $u\in \mathcal{C}([0,1/2])$ and $\Psi \in \mathcal{C}_{b}(L^2(0,1))$
\BD
\E\left[\Psi(Y_{c})|Y_{c}=u \text{ on } [0,1/2]\right] = \E[\Psi(\tilde{B}(c,u))],
\ED
where
\BD
\tilde{B}(c,u):=\left\{\begin{array}{ll}
u(\theta), & \theta \in [0,1/2],\\
u(1/2)+\B_{\theta-1/2} - 12 (1/2-\theta)(\theta-1/2)\left(\int_{0}^{1/2}\!\!\!\!\!\B(r)dr+m(u)-c\right), & \theta \in ]1/2,1].
\end{array}\right.
\ED
Then we have:
\BD
\begin{array}{rcl}
\int_{0}^{1/2}\int_{H}\Psi d\sigma_{r,c}^n\ dr &=& \frac{1}{Z^n_{c}}\int_{H} \int_{0}^{1/2}\E\left[\Psi\times f^n\times \tilde{\gamma}^{,n}\left(\tilde{B}(c,u)\right)\right]\tilde{\gamma}^{n}(u)\mu_{c}(du)\ dr,\\
\\
&=& \frac{1}{Z^n_{c}}\int_{H} \int_{0}^{1/2}\E\left[\Psi\times \tilde{\gamma}^{,n}\left(\tilde{B}(c,u)\right)\right]f^n(u(r))\tilde{\gamma}^{n}(u)\mu_{c}(du)\ dr,
\\
&=& \frac{1}{Z^n_{c}}\int_{H} \int_{0}^{1/2}\E\left[\Psi\times \tilde{\gamma}^{,n}\left(\tilde{B}(c,u)\right)\right]\tilde{\sigma}^n_{r,c}(du).
\end{array}
\ED
Arguing as in the proof of section \ref{S:3.2}, it is easy to conclude that the limit exists, which proves \eqref{Eq:5.28} and \eqref{Eq:5.27}.\\
Recall $\eta$ is the limit of $d\eta^n:=f^n(\hat{X}^n_{c}(t,\theta))dtd\theta-f(\hat{X}_{c}(t,\theta))dtd\theta$. We just proved that for all $\delta\leq s \leq t \leq T$
\BE\label{Eq:5.50}
\begin{array}{rcl}
\E\left[\eta(O_{s,t})\right]&\leq& \liminf_{n\rightarrow+\infty}\E\left[\eta^n(O_{s,t})\right]\\
&\leq&\liminf_{n\rightarrow+\infty}\E\left[\int_{O_{s,t}} f^n(\hat{X}^n_{c}(u,\theta))du\ d\theta\right]\\
\\
&\leq&\liminf_{n\rightarrow+\infty}\int_{H}\int_{O_{s,t}} f^n(x(\theta))du\  d\theta\gamma^n(x)\mu_{c}(dx)\\
&=&(t-s)\liminf_{n\rightarrow+\infty}\int_{H}\int_{0}^1d\sigma^n_{r,c}\ dr\\
&<&+\infty.
\end{array}
\EE
Thus the total mass of $O_{s,t}$ for the reflection measure $\eta$ has a finite expectation.\EProof

\subsection{Reflection and Revuz measures}
\begin{Th}\label{Th:5.4}
For all $c>0$, for all $x \in K\cap H_{c}$:\\
i)\ For $\alpha \geq 3$, the reflection measure $\eta^x$ of the strong solution $\left((X(t,x))_{t\geq0},\eta^x,W\right)$ vanishes.\\
ii)\ For $\alpha < 3$, the reflection measure $\eta^x$ of the strong solution $\left((X(t,x))_{t\geq0},\eta^x,W\right)$ does not vanishes.
\end{Th}
\Proof Let $c>0$, $ x\in K \cap H_{c}$, and $\alpha\geq3$. We take the expectation of equation \eqref{Eq:0.10} for the stationnary solution. We obtain for all $0<\delta \leq s \leq t \leq T$, for all $h \in D(A^2)$:
\BE\label{Eq:5.60}
\E\left[ \int_{s}^t\langle \hat{X}_{c}(u),A^2h\rangle du+\int_{s}^t\langle Ah(\theta),f(\hat{X}_{c}(u))\rangle du+ \big\langle Ah,\eta\big\rangle_{O_{s,t}}\right]=0.
\EE
Thanks to Proposition \ref{Prop:5.4}, the expectation of each term of \eqref{Eq:5.60} is finite. So let $k \in D(A)$ with $\bar{k}=0$, taking $h \in D(A^2)$ such that $k=Ah$ as a test function in \eqref{Eq:5.60}, we obtain for all $0<\delta \leq s \leq t \leq T$, for all $k \in D(A)$:
\BEAS
\E\left[ \big\langle \Pi k,\eta\big\rangle_{O_{s,t}}\right]&=&-\E\left[ \int_{s}^t\langle \hat{X}_{c}(u),Ak\rangle du+\int_{s}^t\langle \Pi k(\theta),f(\hat{X}_{c}(u))\rangle du\right]\\
&=&(s-t)\ \E\left[ \langle \hat{X}_{c}(0),Ak\rangle+\langle \Pi k(\theta),f(\hat{X}_{c}(0))\rangle \right]\\
&=&(s-t) \int_{H} \Big(\langle x, Ak \rangle +\left\langle f(x),\Pi k\right\rangle\Big)\nu_{c}(dx).
\EEAS
We use \eqref{Eq:5.14} with $\phi=1$, and prove that for all $k\in D(A)$, for all $0< \delta\leq s \leq t \leq T$:
\BE\label{Eq:5.38}
\E\left[\big\langle\Pi k , \eta\big\rangle_{O_{s,t}}\right]=0.
\EE
Now, as in Lemma \ref{Le:3.1}, $\eta \otimes \PP$ can be decomposed as $\eta\otimes\PP = \Gamma \otimes d\theta$, where $\Gamma$ is a measure on $[0,T]\times \Omega$, so we obtain that for all $0<\delta\leq s \leq t \leq T$, for all $\mathcal{A}\subset \Omega$:
\BD
0 = \E\left[\langle \hat{X}_{c},\eta\rangle_{O_{s,t}}\1_{\mathcal{A}}\right] =\int_{\Omega}\left( \int_{s}^t\Big(\int_{0}^1\hat{X}_{c}(u)(\theta)d\theta\Big)\1_{\mathcal{A}}\right)d\Gamma(u,.)= c\times \Gamma([s,t],\mathcal{A}).
\ED
Since $c>0$, we conclude that for all $0<\delta\leq s \leq t \leq T$, for all $\mathcal{A}\subset \Omega$, $\Gamma([s,t],\mathcal{A})=0$. Thus $\eta\otimes\PP$ is the null measure. Since $\eta$ is a positive measure, we obtain that $\eta$ is the null measure almost surely. Since the law of $\left(X(t,x)_{t\geq 0},\eta^x\right)$ is a regular conditional distribution of the law of $\left(\hat{X}_c,\eta\right)$ given $\hat{X}_{c}(0)=x \in K \cap H_c$, we have proved i) in the Theorem \ref{Th:5.4}.\\
We consider now the logarithmic case and the negative $\alpha$-power case for $\alpha<3$.
\begin{Prop}
The process $\{X_{c}(t,x) , t \geq 0, x \in H_{c}\cap K\}$ is a continuous Hunt process on $K$ with infinite life-time and strong Markov, properly associated with the Dirichlet Form $\mathcal{E}^c$. In particular, $\mathcal{E}^c$ is quasi-regular.
\end{Prop}
The last assertion is a consequence of Theorem IV.5.1 in \cite{MR1214375}, which describes the necessity of quasi regularity of a Dirichlet Form associated with a Markov process.
We now recall the definitions of the theory of Additive Functionals of a Markov process (see \cite{MR1303354}). Consider $\{Y_{c}(t,x) , t \geq 0, x \in H_{c}\cap K\}$ a Hunt process with infinite life-time and strong Markov, properly associated with the Dirichlet Form $\mathcal{E}^c$. We first describe the minimum admissible filtration and the minimum completed admissible filtration. We set:
\BD
\begin{array}{lcl}
\mathcal{F}^0_{\infty}&=& \sigma\{Y_{c}(s),s\in [0,+\infty)\},\\
\\
\mathcal{F}^0_{t}&=& \sigma\{Y_{c}(s),0\leq s\leq t\},\quad \text{ for all } 0\leq t \leq +\infty.
\end{array}
\ED
These filtrations are called the minimum admissible filtrations. We define $E:=\mathcal{C}([0,T],H_{c})$, and denote by $\PP_{x}$ the law of $t\mapsto Y_{c}(t,x)$ on the filtered space $(E,\mathcal{F}^0_{\infty})$ for all $x\in H_{c} \cap K$. We also define $\mathcal{P}$ the set of all probability measures on $H_{c} \cap K \cap \mathcal{C}([0,1])$. For all $m \in \mathcal{P}$, we define the probability measure:
\BD
\begin{array}{rcl}
\PP_{m} : \mathcal{F}^0_{\infty} & \rightarrow & \R\\
\Lambda &\mapsto & \int_{H_{c} \cap K\cap  \mathcal{C}([0,1])} \PP_{x}(\Lambda)m(dx).
\end{array}
\ED
We then denote by $\mathcal{F}_{\infty}^m$ (respectively $\mathcal{F}_{t}^m$) the completion of $\mathcal{F}_{\infty}^0$ (respectively the completion of $\mathcal{F}_{t}^0$ in $\mathcal{F}_{\infty}^m$) with respect to $\PP_{m}$. Finally we set
\BD
\mathcal{F}_{\infty}= \bigcap_{m\in \mathcal{P}} \mathcal{F}_{\infty}^m\quad\text{ and } \quad\mathcal{F}_{t}= \bigcap_{m\in \mathcal{P}} \mathcal{F}_{t}^m,\quad \text{ for all } 0\leq t \leq +\infty.
\ED
These filtrations are the minimum completed admissible filtrations. It is now possible to define the Additive Functionals for the Markov process $Y_{c}$. To avoid useless definitions, we just recall the definition of a continuous additive functional (CAF in abreviation) in the strict sense of $Y_{c}$.
\begin{Def}\label{De:3}
A family of real valued functions  $\mathcal{A}:=(\mathcal{A}_{t})_{t\geq0}$ is called a continuous additive functionnal in the strict sense of $Y_{c}$ if it satisfies the following conditions:
\begin{enumerate}
\item[(a)] $\mathcal{A}_{t}$ is $\mathcal{F}_{t}$-adapted for all $t\geq0$.
\item[(b)] There exists $\Lambda \in \mathcal{F}_{\infty}$ with $\PP_{x}(\Lambda)=1$, for all $x \in H_{c} \cap K \cap \mathcal{C}([0,1])$, such that $\theta_{t}(\Lambda) \subset \Lambda$ for all $t \geq 0$, for all $\omega \in \Lambda$, $t\mapsto \mathcal{A}_{t}(\omega)$ is continuous, $\mathcal{A}_{0}(\omega)=0$ and  for all $t,s \geq 0$:
\BD
A_{t+s}(\omega) = A_{s} (\omega) + A_{t} (\theta_{s}\omega), 
\ED
where $(\theta_{s})_{s\geq0}$ is the time-translation semigroup on $E$.
\end{enumerate}
Moreover, by a positive continuous additive functional (PCAF in abreviation) in the strict sense of $Y_{c}$, we mean a CAF in the strict sense of $Y_{c}$ such that:
\begin{enumerate}
\item[(c)] For all $\omega \in \Lambda, t \mapsto \mathcal{A}_{t} (\omega)$ is non-decreasing.
\end{enumerate}
\end{Def}
If $\mathcal{A}$ is a linear combination of PCAFs in the strict sense of $Y_{c}$ , the Revuz-measure 
of $\mathcal{A}$ is a Borel signed measure $m$ on $K$ such that for all $\Phi, \Psi \in \mathcal{C}_{b}(H_{c})$: 
\BD
\int_{H_{c}} \Phi(x) \ \E\left[\int_{0}^{+\infty}\!\!\exp(-t)\Psi(Y_{c}(t,x))d\mathcal{A}_{t}\right]\nu_{c}(dx) = \int_{H_{c}} \E\left[\int_{0}^{+\infty}\!\!\exp(-t)\Phi(Y_{c}(t,x))dt\right]\Psi(x) m(dx).
\ED
Notice that there exists a correspondence between Revuz-measures and PCAF. We refer to Chapter 5 in \cite{MR1303354} and Chapter VI in \cite{MR1214375} for all basic definitions and details. In particular the definition of a martingale additive functional (MAF in abreviation), the notion of the energy of an AF, and the quasi-sets.\\
$X_{c}$ does not satisfy suitable properties to compute Revuz-measures of PCAFs in the strict sense of $X_{c}$. Thus we will use a family of process $(Y^{\delta}_{c})_{\delta>0}$ such that:
\BD
Y^{\delta}_{c}(t,x) = X_{c}(t+\delta,x), \quad \text{ for all } x \in H_{c}\cap K, \text{ for all }t\geq 0, \text{ for all } \delta>0. 
\ED
Set $\delta>0$. Let $k\in D(A^2)$, set $h\in D(A)$ such that $Ak=h$ and set $\mathcal{V}: H_{c} \cap K\mapsto \mathcal{V}(x):=\langle x,k\rangle$. Since the Dirichlet form $(\mathcal{E}_{c},D(\mathcal{E}_{c}))$ is quasi-regular, we can apply the Fukushima decomposition (see Theorem VI.2.5 in \cite{MR1214375}). We state that there exists a MAF of finite energy $M^{[\mathcal{V}]}$ and a CAF of zero energy $N^{[\mathcal{V}]}$ such that for $\mathcal{E}^{c}$-quasi every $x$:
\BE\label{Eq:5.39}
\mathcal{V}(Y^{\delta}_{c}(t,x)) - \mathcal{V}(Y^{\delta}_{c}(0,x))= M^{[\mathcal{V}]}_{t} + N^{[\mathcal{V}]}_{t}, \quad t \geq 0, \ \PP^{\delta}_{x}-a.s,
\EE
with obvious notations for $\PP^{\delta}_{x}$. $M^{[\mathcal{V}]}$ and $N^{[\mathcal{V}]}$ can be extended to CAF and MAF in the strict sense of $X_{c}$, which we still denote $M^{[\mathcal{V}]}$ and $N^{[\mathcal{V}]}$, such that $M^{[\mathcal{V}]}$ is a $\PP_{x}$-martingale and \eqref{Eq:5.39} holds for all $x \in H_{c}\cap K$. We have he following expression:
\BD
\frac{1}{2}\int_{\delta}^t\langle Bh,dW\rangle = M^{[\mathcal{V}]}_{t}(Y^{\delta}_{c}(\cdot,x)), \quad \text{ for all } t \geq 0, x \in K, \text{ almost surely}.
\ED
Moreover $N^{[\mathcal{V}]}$ is a linear combination of PCAFs in the strict sense of $Y^{\delta}_{c}$ such that for all $t \geq 0$, $x \in K$, almost surely:
\BD
\frac{1}{2}\int_{\delta}^t\left(\langle X_{c}(s,x), Ah\rangle + \langle f(X_{c}(s,x),\Pi h\rangle \right)ds +
\frac{1}{2}\int_{0}^{1} \Pi h(\theta)\eta^x([\delta,t],d\theta) = N^{[\mathcal{V}]}_{t}(Y^{\delta}_{c}(\cdot,x))
\ED
and its Revuz measure is:
\BE
\frac{1}{2}\left(\langle z, Ah\rangle + \langle f(z),\Pi h\rangle \right) \nu_{c}(dz) + \frac{1}{2}\int_{0}^1\Pi h(r) dr\ \gamma\ d\Sigma^c_{r}.
\EE
To prove the last assertion, it suffices to remark that for all $\Phi, \Psi \in \mathcal{C}_{b}(H_{c})$: 
\BD
\begin{array}{ll}
&2\int_{H_{c}} \Phi(x)\ \E\left[\int_{0}^{+\infty}\exp(-t)\Psi(Y^{\delta}_{c}(t,x))dN^{[\mathcal{V}]}_{t}\right]\nu_{c}(dx)\\
=&\int_{H_{c}}\Psi(x)\ \E\left[ \int_{0}^{+\infty}\!\!\exp(-t)\Phi(Y^{\delta}_{c}(t))dt \right]\left(\langle x, Ah\rangle + \langle f(x),\Pi h\rangle+ \int_{0}^1\Pi h(r) dr\ \gamma\ d\Sigma^c_{r}\right)\nu_{c}(dx).
\end{array}
\ED
Using the same arguments, we remark that there exists a CAF in the strict sense of $Y^{\delta}_{c}$ whose Revuz-measure is
\BE
\frac{1}{2}\left(\langle z, Ah\rangle + \langle f(z),\Pi h\rangle \right) \nu_{c}(dz).
\EE
Since $X_{c}$ is a solution of the equation \eqref{Eq:0.10} in the sense of Definition \ref{De:2}, we obtain that there exists $\mathcal{A}^{[\mathcal{V}]}$ a linear combination of PCAFs in the strict sense  of $X_{c}$
such that:
\BD
\frac{1}{2}\int_{0}^{1} \Pi h(\theta)\eta^x([\delta,t],d\theta) = \mathcal{A}^{[\mathcal{V}]}_{t}(Y^{\delta}_{c}(\cdot,x)), \quad \text{ for all } t \geq 0, x \in K, \text{ almost surely},
\ED
and its Revuz-measure is:
\BE
\frac{1}{2}\int_{0}^{1} \Pi h(r)  dr\ \gamma\ d\Sigma_{r}^c.
\EE 
Finally, we have the following equality:
\BE
\int_{H} \ \E\left[\int_{0}^{+\infty}\exp(-t)\int_{0}^{1}\Pi h(\theta)\eta^x(\delta+dt,d\theta)\right]\nu_{c}(dx) = \int_{0}^{1}\int_{H} \Pi h(r) \gamma d\Sigma_{r}^c \ dr.
\EE
And the reflection measure $\eta^x$ cannot be identically equal to zero.
\newpage


\bibliographystyle{abbrv}

\bibliography{Bibliludo}
\end{document}